\documentclass[10pt,a4paper]{article}

\usepackage{amsmath,amssymb, bbm}
\usepackage{color}
\newcommand{\C}{\mathbb{C}}
\newcommand{\Ct}{\mathbb{C}_t}
\newcommand{\CT}{\mathbb{C}_T}

\newcommand{\Dtt}{\tilde{\mathbb{D}_t}}

\newcommand{\DTt}{\tilde{\mathbb{D}_T}}
\newcommand{\R}{\mathbb{R}}

\newcommand{\Z}{\mathbb{Z}}

\newcommand{\Bo}{\mathbb{B}_0}
\newcommand{\Co}{\mathbb{C}_0}
\newcommand{\Do}{\mathbb{D}_0}

\def\cA{{\mathcal A}}
\def\cC{{\mathcal C}}
\def\cB{{\mathcal B}}

\def\cF{{\mathcal F}}
\def\cP{{\mathcal P}}
\def\cQ{{\mathcal Q}}
\newcommand{\ee}{\varepsilon}
\renewcommand{\aa}{\alpha}

\renewcommand{\div}{{\rm div}\,}

\newcommand{\Sum}{\displaystyle \sum}
\newcommand{\Hs}{\dot{H^s}}

\newcommand{\Gjk}{\mathcal{G}_{j,k}}

\def\d{\partial}

\def\ddj{\dot \Delta_j}

\def\ddk{\dot \Delta_k^v}
\def\tilde{\widetilde}
\def\hat{\widehat}
\newcommand{\D}{\Delta}
\newcommand{\DF}{\Delta_F}

\newcommand{\n}{\nabla}
\newcommand{\G}{\Gamma}

\newcommand{\Fe}{F_{ext}}
\newcommand{\Om}{\Omega}
\newcommand{\Ome}{\Omega_\varepsilon}
\newcommand{\tOm}{\tilde{\Omega}_{QG}}

\newcommand{\ve}{v_\ee}
\newcommand{\we}{w_\ee}
\newcommand{\woe}{w_{0,\ee}}
\newcommand{\voe}{v_{0,\ee}}
\newcommand{\Ue}{U_\ee}
\newcommand{\Uqg}{U_{\ee,QG}}
\newcommand{\ub}{\bar{u}}
\newcommand{\Uosc}{U_{\ee, osc}}
\newcommand{\Uoe}{U_{0,\ee}}
\newcommand{\Uoqg}{U_{0,\ee,QG}}
\newcommand{\Uoosc}{U_{0,\ee, osc}}
\newcommand{\tUqg}{\tilde{U}_{QG}}
\newcommand{\tvqg}{\tilde{v}_{QG}}
\newcommand{\Tqg}{\tilde{T}_{QG}^*}
\newcommand{\Te}{T_\ee^*}
\newcommand{\tUoqg}{\tilde{U}_{0, QG}}

\newcommand{\We}{W_\ee}
\newcommand{\Weh}{W_\ee^h}
\newcommand{\Wei}{W_\ee^{inh}}
\newcommand{\de}{\delta_{\ee}}
\newcommand{\Phie}{\Phi_\ee}
\newcommand{\Thee}{\theta_\ee}

\renewcommand{\Re}{R_\ee}
\newcommand{\re}{r_\ee}

\newcommand{\cPrr}{\cP_{\re, \Re}}

\newtheorem{thm}{Theorem}
\newtheorem{lem}{Lemma}

\newtheorem{prop}{Proposition}
\newtheorem{defi}{Definition}

\newtheorem{rem}{Remark}

\title{Asymptotics for the rotating fluids and primitive systems with large ill-prepared initial data in critical spaces}

\author{Fr\'ed\'eric Charve\footnote{Univ Paris Est Creteil, CNRS, LAMA, F-94010 Creteil, 2 Univ Gustave Eiffel, LAMA, F-77447 Marne-la-Vall\'ee, France. E-mail: frederic.charve@u-pec.fr}}

\date{}
\begin{document}

\maketitle

\begin{abstract} In this article we study the lifespan and asymptotics (in the large rotation and stratification regime) for the Primitive system for highly ill-prepared initial data in critical spaces. Compared to our previous works, we simplified the proof and made it adaptable to the Rotating fluids system with highly ill-prepared initial data decomposed as a sum of 2D horizontal part and a very large 3D part. We also provide explicit convergence rates.
\end{abstract}

\section{Introduction}

\subsection{Geophysical fluids}

Geophysical fluids dynamics are influenced by two concurrent "forces": the Coriolis force (induced by the rotation of the Earth around its axis), and the vertical stratification of the density (induced by gravity) in a way that can be measured through the Rossby and Froude numbers, namely $Ro$ and $Fr$. The smaller they are, the more influent are these two forces. In this article we will consider, in the whole space, first the Primitive System (sometimes also called Primitive Equations)  and seconds the Rotating fluids system, only considering the rotationnal effects.

Let us first introduce the Primitive system: we are interested by the regime where both phenomena are of the same scale (that is we choose $Ro=\ee$ and $Fr=\ee F$ with $F>0$) and we will call $\ee$ the Rossby number and $F$ the Froude number. The system is written as follows:
\begin{equation}
\begin{cases}
\d_t \Ue +\ve\cdot \n \Ue -L \Ue +\frac{1}{\ee} \cA \Ue=\frac{1}{\ee} (-\n \Phie, 0),\\
\div \ve=0,\\
{\Ue}_{|t=0}=U_{0,\ee}.
\end{cases}
\label{PE}
\tag{$PE_\ee$}
\end{equation}
The unknowns are on one hand $\Ue =(\ve, \Thee)=(\ve^1, \ve^2, \ve^3, \Thee)$, where $\ve$ denotes the velocity of the fluid and $\Thee$ the scalar potential temperature (linked to the density, temperature and salinity), and on the other hand $\Phie$, which is called the geopotential and gathers the pressure term and the centrifugal force. The diffusion operator $L$ is defined by
$$
L\Ue \overset{\mbox{def}}{=} (\nu \D \ve, \nu' \D \Thee),
$$
where $\nu, \nu'>0$ denote the kinematic viscosity and thermal diffusivity (both will be considered as viscosities). The last term $\ee^{-1}\cA$ gathers the rotation and stratification effects and the matrix $\cA$ is defined by
$$
\cA \overset{\mbox{def}}{=}\left(
\begin{array}{llll}
0 & -1 & 0 & 0\\
1 & 0 & 0 & 0\\
0 & 0 & 0 & F^{-1}\\
0 & 0 & -F^{-1} & 0
\end{array}
\right).
$$
The rotating fluids system is what we obtain if we only consider the velocity and neglect the last line and column of $\cA$, it is written as follows:

\begin{equation}
 \begin{cases}
  \d_t \ve +\ve\cdot \n \ve-\nu \D \ve +\frac{e_3\wedge \ve}{\ee} =-\n p_\ee,\\
  \div \ve=0,\\
  {\ve}_{|t=0}= v_0.
 \end{cases}
\label{RF} \tag{$RF_\ee$}
\end{equation}

Both system are variations of the Navier-Stokes system, but each of them features a  special structure brought by their respective limit systems as $\ee$ goes to zero: the QG/oscillating structure for \eqref{PE}, and the 2D-3D structure for \eqref{RF}. More details are given in the following parts.
\\
We will use the same notations as in \cite{FCPAA, FCcompl}: for $s\in \R$ and $T>0$ we define the space:
$$
\dot{E}_T^s=\mathcal{C}_T(\Hs (\R^3)) \cap L_T^2(\dot{H}^{s+1}(\R^3)),
$$
endowed with the following norm (For \eqref{PE} $\nu_0=\min(\nu,\nu')$, and for \eqref{RF} $\nu_0=\nu$):
$$
\|f\|_{\dot{E}_T^s}^2 \overset{def}{=}\|f\|_{L_T^\infty \Hs }^2+\nu_0 \int_0^T \|f(\tau)\|_{\dot{H}^{s+1}}^2 d\tau,
$$
where $H^s(\R^3)$ and $\dot{H}^s(\R^3)$ respectively denote the inhomogeneous and homogeneous Sobolev spaces of index $s\in \R$.\\
When $T=\infty$ we simply denote $\dot{E}^s$ and the corresponding norm is taken over $\R_+$ in time.

\subsection{Primitive system: strong solutions, limit system and QG/osc decomposition}

As emphasized in \cite{FCPAA, FCcompl}, thanks to the skew-symmetry of $\cA$, the classical energy method used to study the Navier-Stokes system (based on $L^2$ or $H^s/\Hs $ inner products) do not "see" the penalized terms and are easily adapted to System \eqref{PE}. In the present work we will only focus on the strong solutions provided by the Fujita-Kato theorem: for any fixed $\ee>0$, if $\Uoe\in \dot{H}^\frac12$ there exists a unique local-in-time strong solution, $\Ue$, defined on $[0,\Te[$ and such that for any $T<\Te$, $\Ue\in \dot{E}_T^\frac12$. Note that the solution is global (that is $\Te=+\infty$) when the initial norm $\|\Uoe\|_{\dot{H}^\frac12}$ is bounded by $c\nu_0$ for some small $c>0$. Let us also recall the following blow up criterion: if the lifespan $\Te$ is finite then:
\begin{equation}
 \int_0^{\Te} \|\n \Ue(\tau)\|_{\dot{H}^\frac{1}{2}(\R^3)}^2 d\tau=\infty.
 \label{critereexpl}
\end{equation}
Moreover, if in addition $\Uoe \in \Hs$ (for some fixed $s\in]-\frac32, \frac32[$) then we can propagate the regularity as done for the Navier-Stokes system: $\Ue \in E_T^s$ for any $T<\Te$. All these results are true wether $F=1$ (non-dispersive regime, we refer to \cite{Chemin2, FCF1}) or $F\neq 1$ (dispersive regime).

In the present work, for a fixed $F\neq 1$, our interest is to study the convergence (and obtain convergence rates) when $\ee$ goes to zero (that is for fast rotating and highly stratified systems) in the case of \emph{ill-posed large initial data} (see below for more details).

We refer to \cite{Chemin2, FC, FC2, FCPAA, FCcompl} for studies of the limit system as the small parameter $\ee$ goes to zero and we will only recall here that this limit system is a transport-diffusion system coupled with a Biot-Savart inversion law and is called the 3D Quasi-geostrophic system:
\begin{equation}
\begin{cases}
\d_t \tOm +\tvqg .\n \tOm -\G \tOm =0,\\
\tUqg =(\tvqg ,\tilde{\theta}_{QG})=(-\partial_2, \partial_1, 0, -F\partial_3) \DF^{-1} \tOm,
\end{cases}
\label{QG}\tag{$QG$}
\end{equation}
where we set $\DF=\d_1^2 +\d_2^2 +F^2 \d_3^2$, and the operator $\G$ is defined by:
$$
\G \overset{def}{=} \D \DF^{-1} (\nu \d_1^2 +\nu \d_2^2+ \nu' F^2 \d_3^2),
$$
The quantity $\tOm=\d_1 \tvqg ^2 -\d_2 \tvqg ^1 -F \d_3 \tilde{\theta}_{QG}$ is called the potential vorticity and led by this  limit system we introduce the following decomposition. Let $U=(v, \theta)$ be a 4-dimensional vectorfield, we first define its potential vorticity $\Om(U)$:
$$
\Om(U)\overset{def}{=} \d_1 v^2 -\d_2 v^1 -F\d_3 \theta,
$$
then its orthogonal decomposition into its quasi-geostrophic and oscillating (or oscillatory) parts (in the same spirit as the Leray or Helmholtz decompositions):
\begin{equation}
U_{QG}=\cQ (U) \overset{def}{=} \left(
\begin{array}{c}
-\d_2\\
\d_1\\
0\\
-F\d_3
\end{array}
\right) \D_F^{-1} \Om (U), \quad \mbox{and} \quad U_{osc}=\cP (U) \overset{def}{=} U-U_{QG}.
\end{equation}
\begin{defi}
 \sl{We will say that a vectorfield with four components $U$ is quasi-geostrophic when $U=\cQ U$, and oscillating (or oscillatory) when $U=\cP U$.}
\end{defi}
We refer to  \cite{Chemin2, FC, FC2, FCpochesLp, FCF1, FCPAA, FCcompl}) for more properties of the associated orthogonal projectors $\cQ$ and $\cP$. In particular System \eqref{QG} can be rewritten into:
\begin{equation}
\begin{cases}
\d_t \tUqg +\cQ(\tvqg .\n \tUqg) -\G \tUqg =0,\\
\tUqg =\mathcal{Q} (\tUqg ),\\
{\tilde{U}_{QG|t=0}= \tUoqg.}
\end{cases}
\label{QG1}\tag{$QG$}
\end{equation}

Not only can we adapt the Leray and Fujita-Kato theorems to System \eqref{QG}, but this system also enjoys more "2D"-features as described in Theorem 14 from \cite{FCPAA} (see also \cite{FC2, FC4, FCPAA}): if we make additional low-frequency assumptions, namely $\tUoqg \in H^{\frac12+\delta}$, we obtain global existence in $\dot{E}^0 \cap \dot{E}^{\frac12+\delta}$ (see below for this notation) without any smallness condition on the initial data. In the present article we \emph{only assume that $\tUoqg \in \dot{H}^\frac12 \cap \dot{H}^{\frac12+\delta}$}, and in this case we only rely on the classical Fujita-Kato theorem that we state here in a form including the regularity propagation property (and without assumptions on $\nu,\nu'>0$):

\begin{thm}
 \sl{For any $\tUoqg \in \dot{H}^\frac12 (\R^3)$, there exists a maximal lifespan $\Tqg>0$ and a unique solution $\tUqg \in \dot{E}_t^\frac12$ for all $t<\Tqg$. Moreover
 \begin{itemize}
  \item There exists $c_0>0$ such that if  $\|\tUoqg\|_{\dot{H}^\frac12} \leq c_0 \min(\nu, \nu')$ then $\Tqg=+\infty$ and for any $t\geq 0$,
$$
\|\tUqg(t)\|_{\dot{H}^\frac12}^2+\min(\nu,\nu') \int_0^t \|\tUqg(\tau)\|_{\dot{H}^{\frac32}}^2 d\tau \leq \|\tUoqg\|_{\dot{H}^\frac12}^2.
$$
\item We have the following blow-up criterion:
$$
\int_0^{\Tqg} \|\tUqg(\tau)\|_{\dot{H}^{\frac32}}^2 d\tau <+\infty \Longrightarrow \Tqg=+\infty.
$$
\item Finally, if for some $s\in ]-\frac32, \frac32[$ we have $\tUoqg\in \dot{H}^\frac12 \cap \dot{H}^s$ then for any $t<\Tqg$, $\tUqg \in \dot{E}_t^\frac12 \cap \dot{E}_t^s$ and there exists a constant $C=C_s>0$ such that,
$$
\|\tUqg(t)\|_{\dot{H}^s}^2+\min(\nu,\nu') \int_0^t \|\tUqg(\tau)\|_{\dot{H}^{s+1}}^2 d\tau \leq \|\tUoqg\|_{\dot{H}^s}^2 e^{\frac{C}{\min(\nu,\nu')} \int_0^t \|\tUqg(\tau)\|_{\dot{H}^{\frac32}}^2 d\tau}.
$$
 \end{itemize}
 \label{thQG}
 }
\end{thm}
Going back to System \eqref{PE}, we introduce $\Ome=\Om(\Ue)$, and the usual procedure is then to study separately $\Uqg=\cQ(\Ue)$ and $\Uosc=\cP(\Ue)$. We also decompose the initial data into its oscillating and quasi-geostrophic parts: $\Uoe =\Uoosc +\Uoqg$. We will assume that the $QG$-part converges to some quasi-geostrophic vectorfield $\tUoqg$, and that the oscillating part is very large (in terms of the Rossby number $\ee$) we say that such an initial data is \textbf{ill-prepared}. We refer to \cite{FCcompl} for a small survey concerning recent results about this system and to \cite{Chemin2, BMN5, FC, FCpochesLp, FCF1, FCPAA} for more details.

In the present article, we ask from now on that $F\neq 1$ and mainly focus on the case $\nu=\nu'$ (we will sometimes make remarks about the results in the case $\nu\neq \nu'$) to extend our results for less regular initial data. More precisely, in the continuity of \cite{FCPAA, FCcompl}, we are interested in showing that for very large ill-posed, and less regular, initial data, we are still able to show that the solutions of \eqref{PE} converge to the solution of \eqref{QG}, and provide a convergence rate according to the following sketchy statement:
\begin{thm} \textbf{(Rough statement of the results)}
 \sl{For large ill-posed initial data in $\dot{H}^\frac12 \cap \dot{H}^{\frac12+\delta}$ (initial oscillating part of size $\ee^{-\gamma}$) the lifespan $\Te$ can be made as close to $\Tqg$ as desired provided that the Rossby number $\ee$ is small enough. Moreover, we show that $\||D|^\beta(\Ue-\tUqg)\|_{L^2 L^\infty}$ is of size $\ee^{\beta'}$ for some small $\beta, \beta'>0$. We can reach $\beta=0$ with additional low frequency assumptions on the initial data.}
\end{thm}
We also simplified the proofs so that we can adapt them to prove similar results for the case of the Rotating fluids system.

\subsection{Primitive system: auxiliary systems and statement of the results}
\label{IntroPE}
As in \cite{FCPAA, FCcompl}, we will not be able to estimate directly $\Ue-\tUqg$ and will need to introduce auxiliary systems that will also help us stating our results. With the usual notation, for $f:\R^3 \rightarrow \R^4$, $f\cdot \n f=\sum_{i=1}^3 f^i \d_i f$, let us first rewrite \eqref{QG} as follows:
\begin{equation}
\begin{cases}
{\d_t \tUqg -\G \D \Ue+\frac{1}{\ee} \mathbb{P} \mathcal{A} \tUqg=-\mathbb{P}(\tUqg.\n \tUqg) +G,}\\
{\tilde{U}_{QG|t=0}= \tUoqg.}
\end{cases}
\label{QG3}\tag{$QG$}
\end{equation}
where $G$ is the following divergence-free and potential vorticity-free vectorfield defined as
\begin{equation}
G= G^b + G^l \overset{def}{=}
\mathbb{P} \mathcal{P} (\tUqg. \n \tUqg) -F(\nu-\nu')\Delta \Delta_F^{-2}
\operatorname{} \left(
  \begin{array}{c}
-F \partial_2 \partial_3^2 \\
F \partial_1 \partial_3^2 \\
0\\
(\partial_1^2+\partial_2^2)\partial_3
\end{array} \right) \tOm.
\label{G}
\end{equation}
In \cite{FCPAA, FCcompl}, we then considered the solution $\We$ of the following linear system :
\begin{equation}
 \begin{cases}
  \d_t \We -\G \D \We +\frac{1}{\varepsilon} \mathbb{P} \mathcal{A} \We = -G,\\
  {\We}_{|t=0}=\Uoosc.
 \end{cases}
\label{We}
 \end{equation}
In the present paper the fact that we only assume $\tUoqg \in \dot{H}^\frac12 \cap \dot{H}^{\frac12+\delta}$ makes a major difference as the term $G$ is now much less handy to manipulate (not $L^1$ in time anymore, more details below), which suggests to split $\We$ as  follows $\We= \Weh+\Wei$ with:
\begin{equation}
 \begin{cases}
  \d_t \Weh -\G \D \Weh +\frac{1}{\varepsilon} \mathbb{P} \mathcal{A} \Weh = 0,\\
  {\Weh}_{|t=0}=\Uoosc,
 \end{cases}
 \quad \mbox{and} \quad
 \begin{cases}
  \d_t \Wei -\G \D \Wei +\frac{1}{\varepsilon} \mathbb{P} \mathcal{A} \Wei = -G,\\
  {\Wei}_{|t=0}=0,
 \end{cases}
 \label{Whinh}
\end{equation}
each one with a different behaviour and requiring different Strichartz estimates as outlined in Proposition \ref{PropestimG} below.
From now on, we switch to the case $\nu=\nu'$, inducing the following simplifications: the non-local operator $\G$ turns into $\nu \D$ and $G=G^b$ (we refer to \cite{FCPAA, FCcompl} for more details). Next, we define
$$
\de=\Ue -\tUqg -\Weh -\Wei,
$$
and focus on the system it satisfies, that we write here:
\begin{equation}
 \begin{cases}
  \d_t \de -\nu \D \de +\frac{1}{\ee} \mathbb{P} \mathcal{A} \de = \Sum_{i=1}^{10} F_i,\\
  {\de}_{|t=0}= \Uoqg-\tUoqg,
 \end{cases}
\label{GE}
\end{equation}
with:
\begin{equation}
\begin{cases}
 F_1 \overset{def}{=}-\mathbb{P}(\de \cdot \n \de), \qquad F_2 \overset{def}{=}-\mathbb{P}\Big(\de \cdot \n (\tUqg+\Wei)\Big), \qquad F_3 \overset{def}{=}-\mathbb{P}\Big((\tUqg+\Wei) \cdot \n \de\Big),\\
 F_4 \overset{def}{=}-\mathbb{P}(\de \cdot \n \Weh),\qquad F_5 \overset{def}{=}-\mathbb{P}(\Weh \cdot \n \de),\qquad F_6 \overset{def}{=}-\mathbb{P}(\tUqg \cdot \n \Wei),\\
 F_7 \overset{def}{=}-\mathbb{P}\Big((\tUqg+\Wei) \cdot \n \Weh\Big),\qquad F_8 \overset{def}{=}-\mathbb{P}\Big(\Weh \cdot \n (\tUqg+\Wei)\Big)\\
 F_9 \overset{def}{=}-\mathbb{P}\Big(\Wei \cdot \n(\tUqg+\Wei)\Big),\qquad F_{10} \overset{def}{=}-\mathbb{P}(\Weh \cdot \n \Weh).
\end{cases}
 \label{F110}
\end{equation}
Let us state our first result:
\begin{thm} \textbf{(No smallness assumption)}
 \sl{Assume $F\neq 1$ and $\nu=\nu'$ and let $\Tqg$ and $\Te$ be the lifespan of $\tUqg$ and $\Ue$ as introduced previously.
 \begin{enumerate}
  \item For any $T<\Tqg$, $\C_0\geq 1$, $\delta\in]0,\frac16]$ and any $\aa_0>0$, there exist $\ee_T, m_T>0$ and $\DTt\geq 1$ (depending on $F, \nu,\Co, \delta, \alpha_0$ and $T$) such that for all $\ee\in]0,\ee_T]$ and all divergence-free initial data $\Uoe= \Uoqg + \Uoosc \in \dot{H}^{\frac12} \cap \dot{H}^{\frac12 + \delta}$ satisfying the following assumptions:
 \begin{itemize}
  \item $(H_1) \quad$ There exists a quasi-geostrophic vectorfield $\tUoqg\in \dot{H}^{\frac12} \cap \dot{H}^{\frac12 + \delta}$ such that
  $$
   \begin{cases}
    \|\Uoqg-\tUoqg\|_{\dot{H}^{\frac12} \cap \dot{H}^{\frac12 + \delta}}\leq \C_0 \ee^{\aa_0},\\
    \|\tUoqg\|_{\dot{H}^{\frac12} \cap \dot{H}^{\frac12 + \delta}}\leq \C_0.
   \end{cases}
  $$
\item $(H_2)\quad \|\Uoosc\|_{\dot{H}^{\frac12 + \delta}} \leq m(\ee) \ee^{-\frac{\delta}2}$, with $0<m(\ee)\leq m_T$,
 \end{itemize}
we have $\Te>T$ and with $\Weh, \Wei$ and $\de$ defined as previously,
$$
\|\de\|_{\dot{E}_T^\frac12}\leq \DTt \max(\ee^{\aa_0},\ee^{\frac{\delta}2}, m(\ee)).
$$
\item For any $T<\Tqg$, $\C_0\geq 1$, $\delta\in]0,\frac16]$, $\gamma\in]0,\frac{\delta}{2}[$ and any $\aa_0>0$, if $\eta_0=\frac12(1-\frac{2\gamma}{\delta})$ (or equivalently $\gamma=(1-2\eta_0)\frac{\delta}{2}$), there exist $\ee_T>0$ and $\DTt\geq 1$ (depending on $F, \nu,\Co, \delta, \alpha_0, \gamma$, and $T$) such that for all $\ee\in]0,\ee_T]$ and all divergence-free initial data $\Uoe= \Uoqg + \Uoosc \in \dot{H}^{\frac12} \cap \dot{H}^{\frac12 + \delta}$ satisfying $(H_1)$ and:
 \begin{itemize}
  \item $(H_3)\quad \|\Uoosc\|_{\dot{H}^{\frac12 + \delta}} \leq \C_0 \ee^{-\gamma}$,
 \end{itemize}
then the following results are true:
\begin{enumerate}
 \item $\Te>T$ and for all $s\in[\frac12, \frac12+2\eta_0 \delta[$, we have
\begin{equation}
 \|\de\|_{\dot{E}_T^s}\leq \DTt \ee^{\min\big(\aa_0, \frac{\delta}2-\gamma +\frac12 (\frac12-s)\big)}= \DTt \ee^{\min\big(\aa_0, \frac12 (\frac12 +2\eta_0 \delta-s) \big)}.
\label{estimdeltath2}
 \end{equation}
\item Moreover if, in addition, there exists $c \in ]0,1[$ (assumed to be close to $1$) such that
\begin{itemize}
 \item $(H_4)\quad \|\Uoosc\|_{\dot{H}^{\frac12+ c\delta}\cap \dot{H}^{\frac12 + \delta}} \leq \C_0 \ee^{-\gamma}$,
\end{itemize}
then we can get rid of the oscillations: for all $\eta\in]0,2\eta_0[$, for all $\eta'\in]0,\min(\eta, c)[$, if $\ee\in]0,\ee_T]$ ($\ee_T$ and $\DTt$ now also depend on $c,\eta,\eta'$) then:
$$
\big\||D|^{\eta'\delta}(\Ue-\tUqg)\big\|_{L_T^2 L^\infty} \leq \DTt \ee^{\min\big(\aa_0, (\eta_0-\frac{\eta}2)\delta\big)}.
$$
\item Finally, with more low-frequency regularity on the initial oscillating part, that is
\begin{itemize}
 \item $(H_5)$ $\Uoosc,\; \Uoqg,\; \tUoqg\in \dot{H}^{\frac12-\delta} \cap \dot{H}^{\frac12 + \delta}$, $\Uoosc$ satisfies $(H_4)$, and $(H_1)$ is modified as follows:
\end{itemize}
$$
\begin{cases}
\|\Uoqg-\tUoqg\|_{\dot{H}^{\frac12-\delta} \cap \dot{H}^{\frac12 + \delta}}\leq \C_0 \ee^{\aa_0},\\
    \|\tUoqg\|_{\dot{H}^{\frac12-\delta} \cap \dot{H}^{\frac12 + \delta}}\leq \C_0,
\end{cases}
$$
then for any $T<\Tqg$ and any $k\in]0,1[$ (as close to $1$ as we wish), if $\ee\leq \ee_T$ then \eqref{estimdeltath2} can be extended, for all $s\in[\frac12 -\eta \delta, \frac12[$ (with $0<\eta <2\eta_0$) into:
$$
\|\de\|_{\dot{E}_T^s}\leq \DTt \ee^{\min\big(\aa_0, \frac{1}{2-s}k\eta_0 \delta, \frac12(\frac12 +2\eta_0 \delta-s) \big)},
$$
and finally, we have for all $\ee\leq \ee_T$:
$$
\|\Ue-\tUqg\|_{L_T^2 L^\infty} \leq \DTt \ee^{\frac12 \left(\min(\aa_0, \frac23k\eta_0 \delta) +\min(\aa_0, k\eta_0 \delta) \right)},
$$

\end{enumerate}
 \end{enumerate}
\label{ThPE}
}
\end{thm}
Two immediate extensions can be proved:
\begin{thm}\textbf{(Smallness assumption)}
 \sl{If $(H_1)$ is supplemented with $\|\tUoqg\|_{\dot{H}^\frac12} \leq c_0 \nu$, the previous theorem can be expressed as in \cite{FCPAA, FCcompl}, that is the estimates becomes uniform in time, the constants $\DTt$ become universal constants $\Bo$, "$\Te>T$" becomes "$\Te=+\infty$" and there is no mention to some $T<\Tqg$ anymore.}
\label{ThPEsmall}
 \end{thm}
Finally, when $\nu\neq\nu'$ the result can also be generalized: 
 \begin{thm}\textbf{(Extension of Theorem 21 from \cite{FCPAA}, $\nu\neq \nu'$)}
 \sl{Replacing Assumption (1) by $(H_1)$ but keeping (2) for the oscillating part, allows to extend the result when $\nu\neq \nu'$ for both small and large $\|\tUoqg\|_{\dot{H}^\frac12}$.}
 \label{ThPEext}
\end{thm}
\begin{rem}
 \sl{\begin{enumerate}
      \item Proving the last result requires an adaptation of the Strichartz estimates from \cite{FCPAA} similar to what we did in \cite{FCcompl} and in the present paper in order to improve the condition $r>4$ into $r>2$. The low frequency assumption (2) has to be kept because of truncation arguments.
      \item We tried to simplify the statement of the result from Points 2.b and 2.c compared to \cite{FCcompl}.
     \end{enumerate}
}
\end{rem}
\begin{rem}
 \sl{As we explained, Theorems \ref{ThPE} and \ref{ThPEsmall} are in fact valid for any $\delta<\frac14$ (and the last bound becomes $\DTt \ee^{\min(\aa_0, k\eta_0 \delta)}$). We can prove it with the arguments from \cite{FCPAA, FCcompl} featuring non-local 3D-fractional derivation operators that are adapted neither to anisotropic estimates, nor to 2D-3D products involved in the rotating fluids case. This is why we present here a simplified version of the proof, only holding when $\delta\leq \frac16$, but adapted to prove the corresponding results for the rotating fluids system, which is the object of the following part. The bound for $\delta$ in Theorem \ref{ThPEext} is much smaller, as in \cite{FCPAA}.}
\end{rem}

\subsection{Rotating fluids: auxiliary systems and statement of the results}

As outlined in \cite{CDGG, CDGG2, CDGGbook}, if $v_{0}$ belongs to $L^2(\R^3)$ or $\dot{H}^{\frac12}(\R^3)$ the Leray and Fujita-Kato theorems can be easily adapted but these energy methods fitted to the generic Navier-Stokes system do not take advantage of the special 2D structure induced by strong rotation (when the Rossby number $\ee$ is small).

Moreover, when we consider a more physically relevant initial data of the form $v_0(x)=v_0(x_h,x_3)=\ub_0(x_h)+w_0(x)$ (where $x_h=(x_1,x_2)$ denotes the horizontal variable, and both parts have three components and are divergence-free), the previous results have to be adapted and Chemin, Desjardins, Gallagher and Grenier (we refer to \cite{CDGG, CDGGbook}) first introduce, as a candidate for the limit of the solutions of System \eqref{RF} when $\ee$ goes to zero, $(\ub, \bar{p})=(\ub, \bar{p})(x_h)$ solving the following 2D-Navier-Stokes system (but with three components) :
\begin{equation}
 \begin{cases}
  \d_t \ub +\ub\cdot \n \ub -\nu \D \ub=-\n\bar{p},\\
  \div \ub=0,\\
  \ub_{|t=0}= \ub_0.
 \end{cases}
\label{NS2D} \tag{$2D-NS$}
\end{equation}
This system can be rewritten as follows (with $\ub=(\ub_h, \ub_3)$ and the convention that operators acting only on the horizontal variable are written $\n_h$, $\div_h$ and $\D_h$)
\begin{equation}
 \begin{cases}
  \d_t \ub_h +\ub_h\cdot \n_h \ub_h -\nu \D_h \ub_h=-\n_h \bar{p},\\
  \d_t \ub_3 +\ub_h\cdot \n_h \ub_3 -\nu \D_h \ub_3=0,\\
  \div_h \ub_h=0,\\
  \ub_{|t=0}= \ub_0.
 \end{cases}
\tag{$2D-NS$}
\end{equation}
The fact that there are three components does not change the result compared to the classical 2D-Navier-Stokes system, and we refer to \cite{CDGG, CDGGbook}, for the following result:
\begin{thm}
 \sl{Let $\ub_0\in L^2(\R^2)^3$ such that $\div_h \ub_{0,h}=\d_1 \ub_0^1+\d_2 \ub_0^2=0$. There exists a unique global solution $\ub\in \dot{E}^0 (\R^2)^3$. Moreover this solution belongs to $\cC(\R_+, L^2(\R^2))$ and satisfies the equality:
 $$
 \frac12 \|\ub(t)\|_{L^2}^2 +\nu \int_0^t \|\n \ub(\tau)\|_{L^2}^2 d\tau=\frac12 \|\ub_0\|_{L^2}^2.
 $$
 }
 \label{thN2D}
\end{thm}
Then Chemin, Desjardins, Gallagher and Grenier study the following modified Navier-Stokes-type system (formally resulting from considering $\we=\ve-\ub$ and putting the rotation term involving $\ub$ in the pressure gradient):
\begin{equation}
 \begin{cases}
  \d_t \we +\we\cdot \n \we +\we\cdot \n \ub +\ub\cdot \n \we -\nu \D \we +\frac{e_3\wedge \we}{\ee}= -\n q_\ee,\\
  \div \we=0,\\
  {\we}_{|t=0}= w_0.
 \end{cases}
\label{PRF} \tag{$PRF_\ee$}
\end{equation}
Note that as emphasized for System \eqref{PE}, the rotation term disappears when performing any inner product in $L^2$ or a Sobolev space, and the real difference comes here from the additional transport terms which involve products of 2D and 3D functions that require the following Sobolev product laws:
\begin{prop}
 \sl{There exists a constant $C>0$ such that for any $s,t<1$ with $s+t>0$ and any $u\in \dot{H}^s(\R^2)$, $v\in \dot{H}^t(\R^3)$, then $uv\in \dot{H}^{s+t-1}(\R^3)$ and we have:
 $$
 \|uv\|_{\dot{H}^{s+t-1}(\R^3)} \leq C \|u\|_{\dot{H}^s(\R^2)} \|v\|_{\dot{H}^t(\R^3)}.
 $$
 }
 \label{prod2D3D}
\end{prop}
Then they obtain the Leray and Fujita-Kato theorems for a fixed $\ee$:
\begin{thm} (\cite{CDGG, CDGGbook})
 \sl{Let $\ub_0\in L^2(\R^2)^3$ and let $\ub$ be the associated global solution of System \eqref{NS2D}. If $w_0\in L^2(\R^3)$ with $\div w_0=0$ there exists a global weak Leray solution $\we\in \dot{E}^0$ to System \eqref{PRF} satisfying for all $t\geq 0$:
 $$
 \|\we(t)\|_{L^2}^2 +\nu \int_0^t \|\n\we(\tau)\|_{L^2}^2 d\tau \leq \|w_0\|_{L^2}^2 e^{\frac{C}{\nu^2}\|\ub_0\|_{L^2}^2}.
 $$
 Moreover, this solution converges to $0$ (that is $\ve$ converges to $\ub$) in the sense that for any $q\in]2,6[$ and any $T\geq 0$, we have
 $$
 \lim_{\ee \rightarrow 0} \int_0^T \|\we(\tau)\|_{L^q(\R^3)}^2 d\tau=0.
 $$
}
\end{thm}
The product laws also make it possible to adapt the Fujita-Kato theorem to this modified 3D-Navier-Stokes system:
\begin{thm}
 \sl{Under the same notations:
\begin{itemize}
 \item If $w_0\in \dot{H}^\frac12(\R^3)$ with $\div w_0=0$, there exists a unique local strong (Fujita-Kato) solution $\we$ defined on some $[0, \Te[$ and for any $t<\Te$, $\we \in \dot{E}_t^\frac12$.
 \item Moreover we also have the same blow-up criteria as for Navier-Stokes as well as  regularity propagation when in addition $w_0\in \dot{H}^s$ for some $s\in]-\frac32, \frac32[$.
 \item Finally there exists $c_0>0$ and $C=C(\nu,\|\ub\|_{L^2})$ such that if $\|w_0\|_{\dot{H}^\frac12}\leq c_0 \nu$ then $\Te=+\infty$ and $\|\we\|_{\dot{E}^\frac12} \leq C\|w_0\|_{\dot{H}^\frac12}$. 
\end{itemize}
 }
\end{thm}
This allows to construct $\ve=\we+\ub$ that solves \eqref{RF} with the classical Navier-Stokes tools but more can be done when taking advantage of the special features brought by strong rotation and more precisely by the dispersive properties featured by the following system ($\mathbb{P}$ still denotes the classical Leray projector on divergence-free vectorfields):
\begin{equation}
 \begin{cases}
  \d_t \We -\nu \D \We +\frac{1}{\ee}\mathbb{P} (e_3\wedge \We)=0,\\
  {\We}_{|t=0}=w_0.
 \end{cases}
\label{LRF}\tag{$LRF_\ee$}
\end{equation}
The authors prove Strichartz estimates (see Proposition \ref{StriCDGG}) and obtain the following global existence result:
\begin{thm}(\cite{CDGG, CDGGbook})
 \sl{Let $v_0= \ub_0 + w_0$ with $\ub_0 \in (L^2(\R^2))^3$ and $w_0 \in (\dot{H}^\frac12(\R^3))^3$ (both of them divergence-free). There exists $\ee_0>0$ such that for all $\ee \in]0,\ee_0]$, there is a unique global solution $\ve$ to System \eqref{RF} which satisfies (where $\ub$ and $\We$ are the respective unique solutions of \eqref{NS2D} and \eqref{LRF}):
 \begin{itemize}
  \item $\we=\ve-\ub$ solves \eqref{PRF} in the space $\cC_b^0(\R_+, \dot{H}^\frac12) \cap \dot{E}^\frac12$,
  \item $\|\ve-\ub-\We\|_{\dot{E}^\frac12}\underset{\ee \rightarrow 0}{\longrightarrow} 0.$
 \end{itemize}
 }
\end{thm}
\begin{rem}
\sl{\begin{enumerate}
     \item The result does not require any smallness from the initial data (but of course, $\ee_0$ is taylored depending on the size of the initial data.)
     \item In \cite{CDGG2} the authors extend their result in the case of anisotropic viscosity (and possible zero vertical viscosity).
    \end{enumerate}
}
\end{rem}
In the second part of this article, we wish to extend this result in the spirit of what we did with System \eqref{PE}, considering initial data that depend on $\ee$ and are ill-posed in the sense that their norm blow-up when $\ee$ goes to zero. Asking small extra-regularity allows us to prove in this case global existence of solutions and exhibit an explicit convergence rate as a power of the Rossby number. This is the aim of the following result:
\begin{thm}
 \sl{\begin{enumerate}
 \item For any $\C_0\geq 1$, $\delta\in]0,\frac14]$, $c,k\in]0,1[$ (as close as we wish to 1) and $\gamma\in]0,\frac{\delta}{2}[$, if $\eta_0=\frac12(1-\frac{2\gamma}{\delta})$ (put differently $\gamma=(1-2\eta_0)\frac{\delta}{2}$), there exists $\ee_0>0$ and $\Do\geq 1$ (depending on $\nu,\Co, \delta, c,k, \gamma$) such that for all $\ee\in]0,\ee_0]$ and all initial data $v_0= \ub_0 + \woe$ with $\ub_0 \in (L^2(\R^2))^3$ and $\woe \in (\dot{H}^\frac12(\R^3)\cap \dot{H}^{\frac12+\delta}(\R^3))^3$ (both of them divergence-free) satisfying:
\begin{itemize}
  \item $(H_2') \quad \|\woe\|_{\dot{H}^{\frac12+ c\delta}\cap \dot{H}^{\frac12 + \delta}} \leq \C_0 \ee^{-\gamma}$,
 \end{itemize}
then $\Te=+\infty$ and for all $s\in[\frac12, \frac12+2\eta_0 \delta[$ we have:
\begin{equation}
 \|\de\|_{\dot{E}^s}\leq \Do \ee^{k (\frac12 +2\eta_0 \delta-s)}.
 \label{estimdeltaRF}
 \end{equation}
\item Moreover we can get rid of the oscillations: for all $\eta\in]0,2\eta_0[$, $\eta'\in]0,\min(\eta, c)[$, we have for all $\ee\in]0, \ee_0]$ ($\ee_0,\Do$ now also depend on $\eta,\eta'$)
$$
\big\||D|^{\eta'\delta}\we\big\|_{L^2 L^\infty} =\big\||D|^{\eta'\delta}(\ve-\ub)\big\|_{L^2 L^\infty} \leq \Do \ee^{k\delta (\eta_0-\frac12\eta')}.
$$
\item Finally, if we ask more low-frequency regularity on the initial 3D-part, that is $\woe\in \dot{H}^{\frac12-\delta} \cap \dot{H}^{\frac12 + \delta}$ and still satisfies $(H_2')$, then when $s\in[\frac12 -\eta \delta, \frac12[$ (with $0<\eta <2\eta_0 \min\big(1, \frac{1}{k}-1\big)$) \eqref{estimdeltaRF} becomes for all $\ee\in]0, \ee_0]$:
$$
\|\de\|_{\dot{E}^s}\leq \Do \ee^{\frac{k}{2(2-s)} (\frac12 +2\eta_0 \delta-s)},
$$
and we have:
$$
\|\we\|_{L^2 L^\infty} =\|\ve-\ub\|_{L^2 L^\infty} \leq \Do \ee^{\frac56 k\eta_0\delta}.
$$
\end{enumerate}
\label{ThRF}
}
\end{thm}
\begin{rem}
 \sl{
 \begin{enumerate}
  \item Note that Point 2. is slightly better than Point 2.b from Theorem \ref{ThPE}.
  \item Our result also generalizes the works from \cite{HS, IT2, IT5, KLT} as they consider initial data with only 3D part (here the limit is zero, the solution of System \eqref{NS2D} with $\ub_0=0$) and \cite{LT, IMT} which only consider small initial QG-part in the critical space.
  \item Let us mention \cite{Dutrifoy2} devoted to the Euler-Coriolis system, with initial data also decomposed as a sum of a 2D and a 3D functions. 
 \end{enumerate}
 }
\end{rem}

This article will be structured as follows: we begin with energy estimates for $\Wei$ and $\Weh$. Then we focus on the proof of Theorems \ref{ThPE} and \ref{ThRF}. We postponed to the appendix the proofs of the new Strichartz estimates: first the one needed to deal with $\Wei$ and then the anisotropic Strichartz estimates for $\We$. An important feature of the present article is that the proof we present here is much simpler than in \cite{FCPAA, FCcompl} as we do not resort to non-local fractional derivatives operators, but this simpler method is valid for a narrower range for $\delta$ (when $\delta\leq\frac16$ whereas we can reach $\delta<\frac14$ with the arguments from the cited article).

For the sake of conciseness we will only focus on what is new and will often refer to \cite{FC2, FCPAA,FCcompl} about the quasi-geostrophic structure, and to \cite{CDGG, CDGGbook} for the rotating fluids. We also give minimal details about the Littlewood-Paley decomposition and will mostly refer to \cite{Dbook} for an in-depth study.

\section{Proof of Theorem \ref{ThPE}}

\subsection{Estimates on $G^b$ and $\Wei$}

Let us recall that we defined in \eqref{G} the external force term $G$ (which is equal to $G^b$ when $\nu=\nu'$). If $\Wei$ et $\Weh$ are the solutions of the linear systems from \eqref{Whinh}, then $\Weh$ is globally defined, and $\Wei$ is defined on $[0, \Tqg[$.
\begin{prop}
 \sl{Assume that $\tUoqg \in \dot{H}^\frac12 \cap \dot{H}^{\frac12 +\delta}$ (with $\delta>0$).
 \begin{enumerate}
  \item There exists a constant $C>0$ such that the external force term  satisfies for all $t<\Tqg$:
 \begin{equation}
  \begin{cases}
  \vspace{0.3cm}
   \mbox{For all }r>1, \quad \|G^b\|_{L_t^r \dot{H}^{-\frac32+\frac2{r}}} \leq \frac{C}{\nu^\frac1{r}} \|\tUoqg\|_{\dot{H}^\frac12}^2 e^{\frac{C}{\nu} \Ct},\\
   \mbox{For all }r>\frac{1}{1-\frac{\delta}2}, \quad \|G^b\|_{L_t^r \dot{H}^{-\frac32+\frac2{r}+2\delta}} \leq \frac{C}{\nu^\frac1{r}} \|\tUoqg\|_{\dot{H}^{\frac12+\delta}}^2 e^{\frac{C}{\nu} \Ct},\\
  \end{cases}
 \end{equation}
 where $\Ct=\int_0^t \|\tUqg(\tau)\|_{\dot{H}^\frac32}^2 d\tau$.
 \item There exists a constant $C>0$ such that for all $t<\Tqg$ and $s\in[\frac12, \frac12+2\delta]$,
 \begin{equation}
  \|\Wei(t)\|_{\dot{H}^s}^2+\nu \int_0^t \|\Wei(\tau)\|_{\dot{H}^{s+1}}^2 d\tau \leq \frac{C}{\nu^2}\|\tUoqg\|_{\dot{H}^\frac12 \cap \dot{H}^{\frac12+\delta}}^4 e^{\frac{2C}{\nu}\Ct}.
 \end{equation}
 \end{enumerate}
 \label{PropestimG}
 }
\end{prop}
\begin{rem}
 \sl{\begin{enumerate}
      \item It is immediate to prove that for all $t\geq 0$ and  $s\in[\frac12, \frac12+\delta]$,
 \begin{equation}
  \|\Weh(t)\|_{\dot{H}^s}^2+2\nu \int_0^t \|\Weh(\tau)\|_{\dot{H}^{s+1}}^2 d\tau \leq \|\Uoosc\|_{\dot{H}^s}^2.
 \end{equation}
 Except at the end of the bootstrap argument, we will not use these energy estimates for $\Weh$ as only the norm of $\Uoosc$ in $\dot{H}^{\frac12+c\delta}\cap \dot{H}^{\frac12+\delta}$ is controlled, but with a negative power of $\ee$.
      \item On the contrary, we will abundantly use them for $\Wei$ which is a little more regular than $\tUqg$ and everywhere both of these quantities are involded, we will estimate $\Wei$ similarly to $\tUqg$.
      \item If we only control the $\dot{H}^{\frac12+\delta}$-norm of $\Uoosc$, the best we could hope for in terms of uniform in $\ee$ energy estimates for $\Weh$ would be provided by the Strichartz estimates (see the appendix): for all $t\geq 0$ and $\sigma\in ]\frac34 (\frac12+\delta), \frac12+\delta]$,
      $$
      \|\Weh\|_{L_t^\infty L^{\frac6{3-2\sigma}}}^2 +\nu \|\n \Weh\|_{L_t^2 L^{\frac6{3-2\sigma}}}^2 \leq C \nu^{\frac12+\delta-\sigma} \ee^{\frac3{\sigma}(\frac12+\delta-\sigma)} \|\Uoosc\|_{\dot{H}^{\frac12+\delta}}^2.
      $$
      \item In the case of small initial data ($\|\tUoqg\|_{\dot{H}^\frac12}\leq c_0 \nu$) we simply use the bound $\Ct \leq \frac1{\nu}\|\tUoqg\|_{\dot{H}^\frac12} \leq \frac{\Co}{\nu}$.
     \end{enumerate}
 }
\end{rem}
\textbf{Proof: }From the energy estimates given by Theorem \ref{thQG} (as well as the propagation of the $\dot{H}^{\frac12+\delta}$-regularity) we obtain by complex interpolation that for all $t<\Tqg$ and $p\in[2, \infty]$,
\begin{equation}
\begin{cases}
 \|\tUqg\|_{L_t^p \dot{H}^{\frac12+\frac2{p}}} \leq \frac{C}{\nu^\frac1{p}} \|\tUoqg\|_{\dot{H}^\frac12} e^{\frac{C}{\nu}\Ct},\\
 \vspace{0.2cm}
 \mbox{and}\\
 \|\tUqg\|_{L_t^p \dot{H}^{\frac12+\frac2{p}+\delta}} \leq \frac{C}{\nu^\frac1{p}} \|\tUoqg\|_{\dot{H}^{\frac12+\delta}} e^{\frac{C}{\nu}\Ct}.
\end{cases}
\end{equation}
Thanks to the Bernstein Lemma, the paraproduct and remainder laws (we refer to \cite{Dbook} for the Bony decomposition): for any $s_1,s_2\in[\frac12, \frac12+\delta]$, and any $p,q\in[2, \infty]$:
\begin{multline}
 \|G^b\|_{\dot{H}^{s_1+s_2+\frac2{p}+\frac2{q}-\frac52}} \leq C\Big(\|T_{\tUqg} \n \tUqg\|_{\dot{H}^{s_1+s_2+\frac2{p}+\frac2{q}-\frac52}} +\|T_{\n \tUqg} \tUqg\|_{\dot{H}^{s_1+s_2+\frac2{p}+\frac2{q}-\frac52}}\\
 +\Sum_{i,j=1,...,3}\|\div R(\tUqg^i,\tUqg^j)\|_{\dot{B}_{1,2}^{s_1+s_2+\frac2{p}+\frac2{q}-1}}\Big).
\end{multline}
If $p$ satisfies $\frac2{p}<\frac32-s_1$ we can bound the first term as follows:
\begin{multline}
 \|T_{\tUqg} \n \tUqg\|_{\dot{H}^{s_1+s_2+\frac2{p}+\frac2{q}-\frac52}}\\
 \leq C \|\tUqg\|_{\dot{B}_{\infty, \infty}^{s_1+\frac2{p}-\frac32}} \|\n \tUqg\|_{\dot{H}^{s_2+\frac2{q}-1}} \leq C \|\tUqg\|_{\dot{H}^{s_1+\frac2{p}}} \|\tUqg\|_{\dot{H}^{s_2+\frac2{q}}}. 
\end{multline}
When $\frac2{q}<\frac52-s_2$ (which is true as soon as $\delta<1$), the second term satisfies:
\begin{multline}
 \|T_{\n \tUqg} \tUqg\|_{\dot{H}^{s_1+s_2+\frac2{p}+\frac2{q}-\frac52}}\\
 \leq C \|\n \tUqg\|_{\dot{B}_{\infty, \infty}^{s_2+\frac2{q}-\frac52}} \|\tUqg\|_{\dot{H}^{s_1+\frac2{p}}} \leq C \|\tUqg\|_{\dot{H}^{s_1+\frac2{p}}} \|\tUqg\|_{\dot{H}^{s_2+\frac2{q}}}. 
\end{multline}
And as $s_1+s_2+\frac2{p}+\frac2{q}>0$ we easily get that
$$
\Sum_{i,j=1,...,3}|R(\tUqg^i,\tUqg^j)\|_{\dot{B}_{1,2}^{s_1+s_2+\frac2{p}+\frac2{q}}} \leq C \|\tUqg\|_{\dot{H}^{s_1+\frac2{p}}} \|\tUqg\|_{\dot{H}^{s_2+\frac2{q}}}.
$$
To sum up, we just obtained that under the previous notations, if we set $r$ such that $\frac1{r}=\frac1{p}+\frac1{q}$, then $r$ satisfies $\frac2{r}<\frac52-s_1$ (contrary to $p$, $q$ has no constraint) and we have:
\begin{multline}
 \|G^b\|_{L_t^r \dot{H}^{s_1+s_2+\frac2{p}+\frac2{q}-\frac52}} \leq C \|\tUqg\|_{L_t^p \dot{H}^{s_1+\frac2{p}}} \|\tUqg\|_{L_t^q \dot{H}^{s_2+\frac2{q}}}\\
 \leq \frac{C}{\nu^{\frac1{p}+\frac1{q}}} \|\tUoqg\|_{\dot{H}^{s_1}} \|\tUoqg\|_{\dot{H}^{s_2}}e^{\frac{C}{\nu}\Ct}.
 \label{estimGLr}
\end{multline}
Conversely, if $r$ satisfies $\frac2{r}<\frac52-s_1$ can we find $p$ (with $\frac2{p}<\frac32-s_1$) and $q\in[2,\infty]$ such that $\frac1{r}=\frac1{p}+\frac1{q}$? Introducing $\aa\in]0, \frac52-s_1]$ such that $\frac2{r}=\frac52-s_1-\aa$ we would like to simply take $p$ so that $\frac2{p}=\frac32-s_1-\aa$ which is possible if and only if $\aa\in]0, \frac32-s_1]$, so that two cases have to be considered:
\begin{itemize}
 \item If $\aa\in]0, \frac32-s_1[$, setting $q=2$ and $p$ so that $\frac2{p}=\frac32-s_1-\aa$ ensures that $\frac2{p}<\frac32-s_1$,
 \item If $\aa\in[\frac32-s_1, \frac52-s_1[$, then we simply take $q=r$ and $p=\infty$ and the condition on $p$ is once more satisfied.
\end{itemize}
Writing \eqref{estimGLr} when $s_1=s_2=\frac12$ or $\frac12+\delta$ gives the first part of the proposition.
\\

To prove the second point, let us simply perform the innerproduct in $\dot{H}^s$ (for some $s$) of \eqref{Whinh} with $\Wei$: for all $t<\tUqg$,
\begin{multline}
\frac12 \frac{d}{dt}\|\Wei(t)\|_{\dot{H}^s}^2+\nu \|\Wei(t)\|_{\dot{H}^{s+1}}^2\\
\leq \|G^b\|_{\dot{H}^{s-1}} \|\Wei\|_{\dot{H}^{s+1}} \leq \frac{\nu}2 \|\Wei(\tau)\|_{\dot{H}^{s+1}}^2 +\frac{C}{\nu} \|G^b\|_{\dot{H}^{s-1}}^2.
\end{multline}
Notice that due to point 1 (with $r=2$), $s$ can freely live in $[\frac12, \frac12+2\delta]$ and the result easily follows as $\|\tUoqg\|_{\dot{H}^\frac12 \cap \dot{H}^{\frac12+\delta}}= \max\left ( \|\tUoqg\|_{\dot{H}^\frac12}, \|\tUoqg\|_{\dot{H}^{\frac12+\delta}}\right)$. $\blacksquare$

\subsection{Estimates on $\de$}
\label{estimdeltaPE}
We will only focus on Theorem \ref{ThPE} (without smallness assumptions), the proof of Theorem \ref{ThPEsmall} being easier as $\Ct$ is bounded by $\frac{\Co}{\nu}$.

As we outlined in Section \ref{IntroPE}, when $\tUoqg$ is not assumed to be small in $\dot{H}^\frac12$, $\tUqg$ is defined on $[0,\Tqg[$, as well as $\Wei$. Moreover, thanks to the additional regularity assumptions, for all $t<\Tqg$, $\tUqg$ and $\Wei$ belong to $\dot{E}_t^\frac12 \cap \dot{E}_t^{\frac12+\delta}$. Note that $\Ue$ also belongs to the previous space but for $t<\Te$.

Let us fix some $T<\Tqg$, assume that $\ee$ satisfies Conditions \eqref{condeps1} and \eqref{condeps2} (that is $\ee\leq \ee_T$ for some small $\ee_T$) and assume by contradiction that
\begin{equation}
\Te\leq T,
\label{Contrad1}
\end{equation}
then it is finite and in particular by the blow-up criterion \eqref{critereexpl} is true. Now as in \cite{FCPAA, FCcompl} let us define (with $C$ introduced in \eqref{estimfinale})
\begin{equation}
 T_\ee \overset{def}{=} \sup \{t\in[0,\Te[, \quad \forall t'\leq t, \|\de(t')\|_{\dot{H}^\frac12} \leq \frac{\nu}{4C}\},
\end{equation}
If $\ee>0$ is so small that $\|\delta_\ee(0)\|_{\dot{H}\frac12}\leq\Co \ee^{\aa_0} \leq \frac{\nu}{8C}$ then $T_\ee>0$. Now assume by contradiction that:
\begin{equation}
T_\ee<\Te
 \label{Contrad2}.
\end{equation}
Then for all $t\leq T_\ee <\Te \leq T<\Tqg$, performing (for $s\in[\frac12, \frac12+\eta \delta]$) the $\dot{H}^s$-inner product of System \eqref{GE} by $\de$ we have (the external force terms are defined in \eqref{F110}):
$$
\frac12 \frac{d}{dt} \|\de(t)\|_{\dot{H}^s}^2 +\nu \|\n \de(t)\|_{\dot{H}^s}^2 \leq \Sum_{j=1}^{10} (F_j|\de)_{\dot{H}^s}.
$$
Now we bound each term from the r.h.s. The first three ones are treated exactly like in \cite{FCPAA} and there exists a constant $C>0$ such that:
\begin{equation}
 \begin{cases}
 \vspace{0.1cm}
  |(F_1|\de)_{\dot{H}^s}| \leq C \|\de\|_{\dot{H}^\frac12} \|\de\|_{\dot{H}^{s+1}}^2,\\
  \vspace{0.1cm}
  |(F_2|\de)_{\dot{H}^s}| \leq \frac{\nu}{18}\|\de\|_{\dot{H}^{s+1}}^2 +\frac{C}{\nu} \left(\|\tUqg\|_{\dot{H}^\frac32}^2 +\|\Wei\|_{\dot{H}^\frac32}^2\right)\|\de\|_{\dot{H}^s}^2,\\
  |(F_3|\de)_{\dot{H}^s}| \leq \frac{\nu}{18}\|\de\|_{\dot{H}^{s+1}}^2 +\frac{C}{\nu^3} \left(\|\tUqg\|_{\dot{H}^\frac12}^2 \|\tUqg\|_{\dot{H}^\frac32}^2 +\|\Wei\|_{\dot{H}^\frac12}^2 \|\Wei\|_{\dot{H}^\frac32}^2\right)\|\de\|_{\dot{H}^s}^2.
 \end{cases}
\label{EstimF123}
\end{equation}
The other terms will be bounded differently: when $s\in[0,1]$, we have $2s=(1-\theta)s+\theta(s+1)$ with $\theta=s$, and $1=(1-\theta')s+\theta'(s+1)$ with $\theta'=1-s$,
\begin{multline}
 |(F_4|\de)_{\dot{H}^s}| \leq C \|F_4\|_{L^2} \|\de\|_{\dot{H}^{2s}} \leq C \|\de\|_{L^6} \|\n \Weh\|_{L^3} \|\de\|_{\dot{H}^{2s}} \leq C \|\n \Weh\|_{L^3} \|\de\|_{\dot{H}^1} \|\de\|_{\dot{H}^{2s}}\\ \leq C \|\n \Weh\|_{L^3} \|\de\|_{\dot{H}^s} \|\de\|_{\dot{H}^{s+1}} \leq \frac{\nu}{18}\|\de\|_{\dot{H}^{s+1}}^2 +\frac{C}{\nu} \|\n \Weh\|_{L^3}^2 \|\de\|_{\dot{H}^s}^2.
 \label{EstimF4}
\end{multline}
Similarly (using the Sobolev injections and the fact that $\frac32=(1-\theta')s+\theta'(s+1)$ with $\theta'=\frac32-s$ and the Young inequality with $(\frac43, 4)$),
\begin{multline}
 |(F_5|\de)_{\dot{H}^s}| \leq C \|F_5\|_{L^2} \|\de\|_{\dot{H}^{2s}} \leq C \|\Weh\|_{L^6} \|\n \de\|_{\dot{H}^\frac12} \|\de\|_{\dot{H}^{2s}}\\
 \leq C \|\Weh\|_{L^6} \|\de\|_{\dot{H}^s}^\frac12 \|\de\|_{\dot{H}^{s+1}}^\frac32 \leq \frac{\nu}{18}\|\de\|_{\dot{H}^{s+1}}^2 +\frac{C}{\nu^3} \|\Weh\|_{L^6}^4 \|\de\|_{\dot{H}^s}^2.
 \label{EstimF5}
\end{multline}
Next, with the same tools,
\begin{multline}
 |(F_6|\de)_{\dot{H}^s}| \leq C \|F_6\|_{L^2} \|\de\|_{\dot{H}^{2s}} \leq C \|\tUqg\|_{\dot{H}^\frac12}^\frac12 \|\tUqg\|_{\dot{H}^\frac32}^\frac12 \|\n \Wei\|_{L^3} \|\de\|_{\dot{H}^s}^{1-s} \|\de\|_{\dot{H}^{s+1}}^s\\
\leq C \left(\|\tUqg\|_{\dot{H}^\frac12}^\frac12 \|\tUqg\|_{\dot{H}^\frac32}^{s-\frac12} \|\n \Wei\|_{L^3}\right) \left(\|\tUqg\|_{\dot{H}^\frac32}^{1-s}  \|\de\|_{\dot{H}^s}^{1-s}\right) \|\de\|_{\dot{H}^{s+1}}^s,
\end{multline}
and using the Young inequality with the indices $(2, \frac2{1-s}, \frac2{s})$ we obtain that:
\begin{equation}
 |(F_6|\de)_{\dot{H}^s}| \leq \frac{\nu}{18}\|\de\|_{\dot{H}^{s+1}}^2 +\frac{C}{\nu^{\frac{s}{1-s}}} \|\tUqg\|_{\dot{H}^\frac32}^2 \|\de\|_{\dot{H}^s}^2 +C \|\tUqg\|_{\dot{H}^\frac12} \|\tUqg\|_{\dot{H}^\frac32}^{2s-1} \|\n \Wei\|_{L^3}^2.
 \label{EstimF6}
\end{equation}
Similarly, we obtain:
\begin{multline}
 |(F_7|\de)_{\dot{H}^s}| \leq \frac{\nu}{18}\|\de\|_{\dot{H}^{s+1}}^2 +\frac{C}{\nu^{\frac{s}{1-s}}} \left(\|\tUqg\|_{\dot{H}^\frac32}^2 +\|\Wei\|_{\dot{H}^\frac32}^2\right) \|\de\|_{\dot{H}^s}^2\\
 +C \left(\|\tUqg\|_{\dot{H}^\frac12} \|\tUqg\|_{\dot{H}^\frac32}^{2s-1} +\|\Wei\|_{\dot{H}^\frac12} \|\Wei\|_{\dot{H}^\frac32}^{2s-1}\right)\|\n \Weh\|_{L^3}^2.
 \label{EstimF7}
\end{multline}
Considering the following term (part of $F_8$):
\begin{multline}
 |(\Weh\cdot \n \tUqg|\de)_{\dot{H}^s}| \leq \|\Weh\|_{L^6} \|\tUqg\|_{\dot{H}^\frac32} \|\de\|_{\dot{H}^s}^{1-s} \|\de\|_{\dot{H}^{s+1}}^s\\
 \leq \left(\|\Weh\|_{L^6} \|\tUqg\|_{\dot{H}^\frac32}^s\right) \left(\|\tUqg\|_{\dot{H}^\frac32}^{1-s}\|\de\|_{\dot{H}^s}^{1-s}\right) \|\de\|_{\dot{H}^{s+1}}^s,
\end{multline}
and thanks once more to the Young inequality with the indices $(2, \frac2{1-s}, \frac2{s})$, we can estimate $F_8$ and $F_9$ as follows
\begin{multline}
 |(F_8|\de)_{\dot{H}^s}| +|(F_9|\de)_{\dot{H}^s}| \leq \frac{\nu}{18}\|\de\|_{\dot{H}^{s+1}}^2 +\frac{C}{\nu^{\frac{s}{1-s}}} \left(\|\tUqg\|_{\dot{H}^\frac32}^2 +\|\Wei\|_{\dot{H}^\frac32}^2\right) \|\de\|_{\dot{H}^s}^2\\
 +C \left(\|\tUqg\|_{\dot{H}^\frac32}^{2s} +\|\Wei\|_{\dot{H}^\frac32}^{2s}\right) \left( \|\Wei\|_{L^6}^2 +\|\Weh\|_{L^6}^2\right).
 \label{EstimF89}
\end{multline}
The last term is also bounded with similar arguments:
\begin{multline}
 |(F_{10}|\de)_{\dot{H}^s}| \leq \|\n \Weh\|_{L^3} \left(\|\Weh\|_{L^6} \|\de\|_{\dot{H}^s}^{1-s}\right) \|\de\|_{\dot{H}^{s+1}}^s\\
 \leq \frac{\nu}{18}\|\de\|_{\dot{H}^{s+1}}^2 +\frac{C}{\nu^{\frac{s}{1-s}}} \|\Weh\|_{L^6}^{\frac2{1-s}} \|\de\|_{\dot{H}^s}^2 + C \|\n \Weh\|_{L^3}^2.
\label{EstimF10}
\end{multline}
Gathering \eqref{EstimF123}, \eqref{EstimF4}, \eqref{EstimF5}, \eqref{EstimF6}, \eqref{EstimF7}, \eqref{EstimF89} and \eqref{EstimF10} we end up with:
\begin{equation}
  \frac12\frac{d}{dt} \|\de(t)\|_{\dot{H}^s}^2 +\frac{\nu}2 \|\n \de(t)\|_{\dot{H}^s}^2 \leq C \|\de\|_{\dot{H}^\frac12} \|\de\|_{\dot{H}^{s+1}}^2 + \frac{C}{\nu} M_1(t) \|\de\|_{\dot{H}^s}^2 +C M_2(t),
  \label{estimfinale}
\end{equation}
where
\begin{multline}
 M_1(t) \overset{def}{=} \|\tUqg\|_{\dot{H}^\frac32}^2 \left(1+ \frac1{\nu^{\frac{2s-1}{1-s}}} +\frac1{\nu^2}\|\tUqg\|_{\dot{H}^\frac12}^2\right) +\|\Wei\|_{\dot{H}^\frac32}^2 \left(1+ \frac1{\nu^{\frac{2s-1}{1-s}}} +\frac1{\nu^2}\|\Wei\|_{\dot{H}^\frac12}^2\right)\\
 +\|\n \Weh\|_{L^3}^2 + \frac1{\nu^2} \|\Weh\|_{L^6}^4 +\frac1{\nu^{\frac{2s-1}{1-s}}} \|\Weh\|_{L^6}^{\frac2{1-s}},
\end{multline}
and
\begin{multline}
 M_2(t) \overset{def}{=} \left(\|\tUqg\|_{\dot{H}^\frac12} \|\tUqg\|_{\dot{H}^\frac32}^{2s-1} +\|\Wei\|_{\dot{H}^\frac12} \|\Wei\|_{\dot{H}^\frac32}^{2s-1}\right) \left(\|\n \Weh\|_{L^3}^2+\|\n \Wei\|_{L^3}^2\right)\\
  +\left(\|\tUqg\|_{\dot{H}^\frac32}^{2s} +\|\Wei\|_{\dot{H}^\frac32}^{2s}\right) \left( \|\Wei\|_{L^6}^2 +\|\Weh\|_{L^6}^2\right) +\|\n \Weh\|_{L^3}^2.
\end{multline}
So that for any $t\leq T_\ee <\Te \leq T<\Tqg$, thanks to the Gronwall lemma, the H\"older estimate, the estimates from Theorem \ref{thQG}, Point 2 from Proposition \ref{PropestimG} and $(H_1)$, there exists a constant $\Bo$ depending on $\Co,\nu,C,s$ such that:
\begin{multline}
 \|\de(t)\|_{\dot{H}^s}^2 +\frac{\nu}{2} \int_0^t \|\n \de(\tau)\|_{\dot{H}^s}^2 d\tau \leq C \left(\|\de(0)\|_{\dot{H}^s}^2+\int_0^t M_2(\tau) d\tau\right) e^{\frac{C}{\nu} \int_0^t M_1(\tau) d\tau}\\
 \leq C_{\nu,s} \Bigg[\|\Uoqg-\tUoqg\|_{\dot{H}^s}^2 +\left(\|\n \Weh\|_{L_t^{\frac2{\frac32-s}}L^3}^2 +\|\n \Wei\|_{L_t^{\frac2{\frac32-s}}L^3}^2\right)\\
 \times\left(\|\Wei\|_{L_t^\infty \dot{H}^\frac12} \|\Wei\|_{L_t^2 \dot{H}^\frac32}^{2s-1} +\|\tUqg\|_{L_t^\infty \dot{H}^\frac12} \|\tUqg\|_{L_t^2 \dot{H}^\frac32}^{2s-1}\right)\\
 +\left(\|\Wei\|_{L_t^2 \dot{H}^\frac32}^{2s} +\|\tUqg\|_{L_t^2 \dot{H}^\frac32}^{2s}\right)\left(\|\Weh\|_{L_t^{\frac2{1-s}}L^6}^2 +\|\Wei\|_{L_t^{\frac2{1-s}}L^6}^2 \right) +\|\n \Weh\|_{L_t^2 L^3}^2\Bigg]\\
\times \exp \Bigg(C_{\nu,s}\Bigg\{(1+\|\tUqg\|_{L_t^\infty \dot{H}^\frac12}) \|\tUqg\|_{L_t^2 \dot{H}^\frac32}^2 +(1+\|\Wei\|_{L_t^\infty \dot{H}^\frac12}) \|\Wei\|_{L_t^2 \dot{H}^\frac32}^2\\
+\|\n \Weh\|_{L_t^2 L^3}^2 +\|\Weh\|_{L_t^4 L^6}^4 +\|\Weh\|_{L_t^{\frac2{1-s}}L^6}^{\frac2{1-s}} \Bigg\}\Bigg)\\
\leq \Bo \Bigg[ \ee^{2\aa_0} +\left(\|\n \Weh\|_{L_T^{\frac2{\frac32-s}}L^3}^2 +\|\n \Wei\|_{L_T^{\frac2{\frac32-s}}L^3}^2 +\|\Weh\|_{L_T^{\frac2{1-s}}L^6}^2 +\|\Wei\|_{L_T^{\frac2{1-s}}L^6}^2\right) e^{\frac{3C}{\nu} \CT}\\
+\|\n \Weh\|_{L_t^2 L^3}^2 \Bigg] \times \exp{\Bo \left\{e^{\frac{3C}{\nu} \CT} +\|\n \Weh\|_{L_T^2 L^3}^2 +\|\Weh\|_{L_T^4 L^6}^4 +\|\Weh\|_{L_T^{\frac2{1-s}}L^6}^{\frac2{1-s}} \right\}},
\label{estimCas1}
\end{multline}
where we recall that we introduced $\Ct=\int_0^t \|\tUoqg(\tau)\|_{\dot{H}^\frac32}^2 d\tau \leq \CT<\infty$ in Proposition \ref{PropestimG}.

We can bound the various terms from the previous estimates involving $\Weh$ (and $\Wei$) thanks to the Strichartz estimates provided by Proposition \ref{estimStriPE}.

\subsection{End of the bootstrap argument}

Let us first focus on the proof of the second point from Theorem \ref{ThPE}. For all $s\in[\frac12, \frac12+\eta \delta]$, under Assumption $(H_3)$, simplifying \eqref{estimCas1} with Proposition \ref{estimStriPE} leads for all $t\leq T_\ee$ to: 
\begin{multline}
  \|\de(t)\|_{\dot{H}^s}^2 +\frac{\nu}{2} \int_0^t \|\n \de(\tau)\|_{\dot{H}^s}^2 d\tau\\
  \leq \Bo \left(\ee^{2\aa_0} + \ee^{2\eta_0 \delta} +(\ee^{\frac12+2\eta_0 \delta-s}+\ee^{\frac12+\delta-s}) e^{\frac{5C}{\nu}\CT} \right) e^{\Bo\left(e^{\frac{3C}{\nu}\CT} +\ee^{2\eta_0 \delta} +\ee^{\frac1{1-s}(\frac12+2\eta_0\delta-s)}\right)}.
\end{multline}
We recall that $\eta<2\eta_0<1$ so we have $0<(2\eta_0-\eta)\delta \leq \frac12+2\eta_0 \delta -s\leq 2\eta_0 \delta$, and when $\ee>0$ is so small that:
\begin{equation}
 \ee^{2\eta_0 \delta}\leq \frac12\quad \mbox{and} \quad \ee^{\frac1{1-s}(\frac12+2\eta_0\delta-s)}\leq \ee^{2(2\eta_0-\eta)\delta} \leq \frac12,
 \label{condeps1}
\end{equation}
then the previous estimates turns into:
\begin{multline}
  \|\de(t)\|_{\dot{H}^s}^2 +\frac{\nu}{2} \int_0^t \|\n \de(\tau)\|_{\dot{H}^s}^2 d\tau \leq \DTt \left(\ee^{2\aa_0} +\ee^{2\eta_0 \delta} +\ee^{\frac12+2\eta_0 \delta-s} +\ee^{\frac12+\delta-s}\right)\\
  \leq \DTt \ee^{\min(2\aa_0,\frac12+2\eta_0 \delta-s)} \leq \DTt \ee^{\min(2\aa_0, (2\eta_0-\eta) \delta)},
\label{estimCas1sb}
\end{multline}
where we set for some $t$, $\Dtt \overset{def}{=}\Bo e^{\frac{5C}{\nu}\Ct} e^{\Bo (1+e^{\frac{3C}{\nu}\Ct})}$. Finally if $\ee>0$ also satisfies:
\begin{equation}
 \DTt \ee^{\min(2\aa_0, (2\eta_0-\eta) \delta)} \leq \left(\frac{\nu}{8C}\right)^2,
 \label{condeps2}
\end{equation}
then for all $t\leq T_\ee$, taking $s=\frac12$, we have
$$
\|\de(t)\|_{\dot{H}^\frac12} \leq \frac{\nu}{8C},
$$
which contradicts the definition of $T_\ee$, so that \eqref{Contrad2} is false and $T_\ee=\Te$. Thanks to \eqref{estimCas1sb}, Theorem \ref{thQG} and Proposition \ref{PropestimG}, with $s=\frac12$ for all $t <\Te \leq T<\Tqg$ we have:
\begin{multline}
 \int_0^t \|\n \Ue(\tau)\|_{\dot{H}^\frac12}^2 d\tau \leq \int_0^t \|\n \de(\tau)\|_{\dot{H}^\frac12}^2 d\tau +\int_0^t \|\n \Wei(\tau)\|_{\dot{H}^\frac12}^2 d\tau +\int_0^t \|\n \Weh(\tau)\|_{\dot{H}^\frac12}^2 d\tau\\
\leq \frac1{\nu} \left(\DTt \ee^{2\min(\aa_0, \eta_0 \delta)} +\Co^2 e^{\frac{C}{\nu}\CT} +\frac{C\Co^4}{\nu^2} e^{\frac{2C}{\nu}\CT} +\int_0^T \|\Weh(\tau)\|_{\dot{H}^\frac32}^2 d\tau \right) <\infty,
\end{multline}
which contradicts \eqref{critereexpl} so that \eqref{Contrad1} is also false and $\Te >T$ which concludes the proof of Point 2-a.
\\
To prove Point 1, resuming the previous bootstrap argument, for $s=\frac12$ simplifying \eqref{estimCas1} now under assumption $(H_2)$ leads for all $t\leq T_\ee$ to (when $\ee$ is set so small that $m(\ee)\leq 1$):
\begin{multline}
  \|\de(t)\|_{\dot{H}^\frac12}^2 +\frac{\nu}{2} \int_0^t \|\n \de(\tau)\|_{\dot{H}^\frac12}^2 d\tau\\
  \leq \Bo \left(\ee^{2\aa_0} + m(\ee)^2 +(m(\ee)^2+ \ee^\delta)e^{\frac{3C}{\nu}\CT} \right) e^{\Bo\left(e^{\frac{3C}{\nu}\CT} +m(\ee)^2\right)}\\
  \leq \DTt \left(\ee^{2\aa_0} + m(\ee)^2 +\ee^\delta \right)
\label{estimCas1s}
\end{multline}
and the same method leads to the result.

\subsection{Proof of Point 2.b}
\label{Point2b}
This point is close to the corresponding result from \cite{FCcompl}, but there are two differences: first, we chose in the present article to state a little differently the result and will give some details (even if the proof is close to the one in \cite{FCcompl}), seconds the new term $\Wei$ has to be estimated in addition to $\Weh$.
\\

For any $k\in]0,1[$ (close to $1$), any $\eta\in[0,2\eta_0[$ and any $\eta'\in[0,\eta[$, from \eqref{estimCas1sb} with $s\in\{\frac12, \frac12+\eta\delta\}$ we get:
\begin{multline}
 \||D|^{\eta'\delta} \de\|_{L_T^2 L^\infty} \leq \||D|^{\eta'\delta} \de\|_{L_T^2 \dot{B}_{2,1}^\frac32} \leq \||D|^{\eta'\delta} \de\|_{L_T^2 \dot{H}^{\frac32-\eta'\delta}}^{1-\frac{\eta'}{\eta}} \||D|^{\eta'\delta} \de\|_{L_T^2 \dot{H}^{\frac32+(\eta-\eta')\delta}}^{\frac{\eta'}{\eta}}\\
 \leq \|\de\|_{L_T^2 \dot{H}^\frac32}^{1-\frac{\eta'}{\eta}} \|\de\|_{L_T^2 \dot{H}^{\frac32+\eta\delta}}^{\frac{\eta'}{\eta}} \leq \DTt \ee^{\left((1-\frac{\eta'}{\eta})\min(\aa_0, \eta_0\delta)  +\frac{\eta'}{\eta}\min(\aa_0, (\eta_0-\frac{\eta}2)\delta)\right)}\\
 \leq\DTt \ee^{\min\left(\aa_0, (\eta_0-\frac{\eta}2)\delta\right)}
 \label{Pt2bA}
\end{multline}
Thanks to Proposition \ref{estimStriPE} with $(d,p,r,q)=(\eta'\delta, 2, \infty, 1)$ implies that for $\theta\in[0,1]$,
\begin{equation}
\||D|^{\eta'\delta} \Weh\|_{L^2 L^\infty} \leq \frac{C_{F,\theta}}{\nu^{\frac{1-\theta}4}} \ee^{\frac{\theta}4} \|\Uoosc\|_{\dot{B}_{2,1}^{\frac12+\eta'\delta+\frac{\theta}2}}.
\label{ChoixPE1}
 \end{equation}
Thanks to Lemma \ref{majBs21} with $(\aa,\beta)=(a\frac{\theta}2, b\frac{\theta}2)$ (with $a,b>0$) we can write:
$$
 \|\Uoosc\|_{\dot{B}_{2,1}^{\frac12+\eta'\delta+\frac{\theta}2}} \leq C_{a,b,\theta} \|\Uoosc\|_{\dot{H}^{\frac12+\eta'\delta+\frac{\theta}2(1-a)}}^\frac{b}{a+b} \|\Uoosc\|_{\dot{H}^{\frac12+\eta'\delta+\frac{\theta}2(1+b)}}^\frac{a}{a+b}.
$$
We then choose $a,b>0$ so that
$$
\begin{cases}
 \frac{\theta}2(1-a)=(c-\eta')\delta,\\
 \frac{\theta}2(1+b)=(1-\eta')\delta.
\end{cases}
$$
Take some $b>0$ (to be fixed later), and choose $\theta=\frac2{1+b}(1-\eta')\delta$ then the existence of some $a\in]0,1[$ satisfying the other condition is equivalent to the fact that $b>0$ is so small that $(c-\eta')(1+b) <1-\eta'$ and in that case:
$$
a=1-(1+b)\frac{c-\eta'}{1-\eta'}.
$$
The condition $\theta\leq 1$ is equivalent to $\delta\leq\frac12 \frac{1+b}{1-\eta'}$, which is true when $\delta\leq \frac16$. On the other hand the condition on the "p-index" from Proposition \ref{estimdispRF} is satisfied if and only if $\delta\leq \frac{1+b}{1-\eta'}$ which is implied by the previous condition. Then \eqref{ChoixPE1} turns into
$$
\||D|^{\eta'\delta} \Weh\|_{L^2 L^\infty} \leq C_{F,\Co, \nu, \delta, b,\eta, \eta'} \ee^{\frac{\delta}{2(1+b)}(1-\eta')-\gamma}.
$$
The exponent of $\ee$ also writes $\frac{\delta}{2(1+b)}\left(2\eta_0-\eta'-b(1-2\eta_0)\right)$, which goes to $\delta(\eta_0-\frac{\eta'}2)>\delta(\eta_0-\frac{\eta}2)$ when $b$ goes to zero, so choosing $b>0$ so small that 
$$
\frac{\delta}{2(1+b)}\left(2\eta_0-\eta'-b(1-2\eta_0)\right)= \delta(\eta_0-\frac{\eta}2),
$$
that is $b=\frac{\eta-\eta'}{1-\eta}$, we finally get:
\begin{equation}
\||D|^{\eta'\delta} \Weh\|_{L^2 L^\infty} \leq \Bo \ee^{\delta (\eta_0-\frac12\eta)}.
\label{Pt2bB}
\end{equation}
Similarly, we get that 
\begin{equation}
\big\||D|^{\eta'\delta}\Wei \big\|_{L_T^2 L^\infty} \leq \DTt \ee^{\delta (\eta_0-\frac12\eta)+\gamma}.
\label{Pt2bC}
\end{equation}
Gathering \eqref{Pt2bA}, \eqref{Pt2bB} and \eqref{Pt2bC} ends the proof. $\blacksquare$

\subsection{Proof of Point 2.c}
\label{Point2c}
First let us emphasize that in Section \ref{estimdeltaPE}, two terms have to be estimated differently when $s\in[\frac12-\eta \delta, \frac12]$, namely $F_6$ and $F_7$, because now the exponent satisfies $2s-1<0$ which makes useless \eqref{EstimF6} and \eqref{EstimF7}. Thus we estimate these terms as follows:
\begin{multline}
 |(F_6|\de)_{\dot{H}^s}| \leq C \|F_6\|_{L^2} \|\de\|_{\dot{H}^{2s}} \leq C \|\tUqg\|_{L^\frac6{3-2s}} \|\n \Wei\|_{L^\frac3{s}} \|\de\|_{\dot{H}^s}^{1-s} \|\de\|_{\dot{H}^{s+1}}^s\\
\leq C \left(\|\tUqg\|_{\dot{H}^s} \|\n \Wei\|_{L^\frac3{s}}^\frac1{2-s}\right) \left(\|\n \Wei\|_{L^\frac3{s}}^\frac{1-s}{2-s}  \|\de\|_{\dot{H}^s}^{1-s}\right) \|\de\|_{\dot{H}^{s+1}}^s.
\end{multline}
By the Young inequality with the indices $(2, \frac2{1-s}, \frac2{s})$, we get
\begin{equation}
 |(F_6|\de)_{\dot{H}^s}| \leq \frac{\nu}{18}\|\de\|_{\dot{H}^{s+1}}^2 +\frac{C}{\nu^{\frac{s}{1-s}}} \|\n \Wei\|_{L^\frac3{s}}^\frac2{2-s} \|\de\|_{\dot{H}^s}^2 + \|\tUqg\|_{\dot{H}^s}^2 \|\n \Wei\|_{L^\frac3{s}}^\frac2{2-s}.
 \label{EstimF6bis}
\end{equation}
Similarly we obtain:
\begin{multline}
 |(F_7|\de)_{\dot{H}^s}| \leq \frac{\nu}{18}\|\de\|_{\dot{H}^{s+1}}^2 +\frac{C}{\nu^{\frac{s}{1-s}}} \|\n \Weh\|_{L^\frac3{s}}^\frac2{2-s} \|\de\|_{\dot{H}^s}^2\\
 +\left(\|\tUqg\|_{\dot{H}^s}^2 +\|\Wei\|_{\dot{H}^s}^2\right)\|\n \Weh\|_{L^\frac3{s}}^\frac2{2-s}.
 \label{EstimF7bis}
\end{multline}
so that when $s\in[\frac12-\eta \delta, \frac12]$, the previous functions $M_1$ and $M_2$ are modified according to:
\begin{multline}
 M_1(t) \overset{def}{=} \|\tUqg\|_{\dot{H}^\frac32}^2 \left(1+ \frac1{\nu^{\frac{2s-1}{1-s}}} +\frac1{\nu^2}\|\tUqg\|_{\dot{H}^\frac12}^2\right) +\|\Wei\|_{\dot{H}^\frac32}^2 \left(1+ \frac1{\nu^{\frac{2s-1}{1-s}}} +\frac1{\nu^2}\|\Wei\|_{\dot{H}^\frac12}^2\right)\\
 +\|\n \Weh\|_{L^3}^2 + \frac1{\nu^2} \|\Weh\|_{L^6}^4 +\frac1{\nu^{\frac{2s-1}{1-s}}} \left(\|\Weh\|_{L^6}^{\frac2{1-s}} +\|\n \Weh\|_{L^\frac3{s}}^{\frac2{2-s}} +\|\n \Wei\|_{L^\frac3{s}}^{\frac2{2-s}}\right),
\end{multline}
and
\begin{multline}
 M_2(t) \overset{def}{=} \left(\|\tUqg\|_{\dot{H}^s}^2 +\|\Wei\|_{\dot{H}^s}^2\right) \left(\|\n \Weh\|_{L^\frac3{s}}^{\frac2{2-s}} +\|\n \Wei\|_{L^\frac3{s}}^{\frac2{2-s}}\right)\\
  +\left(\|\tUqg\|_{\dot{H}^\frac32}^{2s} +\|\Wei\|_{\dot{H}^\frac32}^{2s}\right) \left( \|\Wei\|_{L^6}^2 +\|\Weh\|_{L^6}^2\right) +\|\n \Weh\|_{L^3}^2.
\end{multline}
\begin{rem}
 \sl{Note that estimates \eqref{EstimF6bis} and \eqref{EstimF7bis} would be useless in Section \ref{estimdeltaPE} because $\frac2{2-s}\geq 2$ if and only if $s\in[1,2[$, so when $s\in[\frac12, \frac12+\eta \delta] \subset[0,1[$ we cannot use neither Proposition \ref{Propermut}, nor Proposition \ref{Estimdispnu2} which requires $p\geq 2$. We could use the first point of Proposition \ref{injectionLr} with $p=1$ but it would require the use of Lemma \ref{majBs21} which is not possible under Assumption $(H_3)$ alone.}
\end{rem}
As explained in the previous Remark, we are forced to use Proposition \ref{Estimdispnu} with $(d,p,r,q)=(1,\frac2{2-s}, \frac3{s},1)$ and thanks to Assumption $(H_4)$ we will be able to take advantage of Lemma \ref{majBs21} and Proposition \ref{injectionLr}. Thanks to the last part of Proposition \ref{estimStriPE}, when $k,k',k''<1$ are fixed, \eqref{estimCas1} turns, for any $s\in[\frac12-\eta\delta, \frac12]$, into (when $\ee\leq \ee_T$):
\begin{multline}
  \|\de(t)\|_{\dot{H}^s}^2 +\frac{\nu}{2} \int_0^t \|\n \de(\tau)\|_{\dot{H}^s}^2 d\tau\\
  \leq \Bo e^{\frac{4C}{\nu}\CT} \left(\ee^{2\aa_0} + \ee^{\frac2{2-s}k'\eta_0 \delta} +\ee^{\frac12+2\eta_0 \delta-s} \right) e^{\Bo e^{\frac{4C}{\nu}\CT}(1+\ee^{\frac2{2-s}k'\eta_0 \delta})}\\
  \leq \DTt \left(\ee^{2\aa_0} + \ee^{\frac2{2-s}k'\eta_0 \delta} +\ee^{\frac12+2\eta_0 \delta-s} \right),
\end{multline}
and for $s=\frac12-\eta\delta$, we get that
$$
 \|\de\|_{L_T^2 \dot{H}^{\frac32-\eta\delta}}^2 \leq \DTt \left(\ee^{2\aa_0} + \ee^{\frac2{\frac32+\eta\delta}k'\eta_0 \delta} +\ee^{(2\eta_0+\eta) \delta} \right).
$$
If we have chosen $\eta>0$ so small that $\frac2{\frac32+\eta\delta}\geq\frac43k''$, we get that:
$$
 \|\de\|_{L_T^2 \dot{H}^{\frac32-\eta\delta}} \leq \DTt \ee^{\min(\aa_0, \frac23k'k''\eta_0 \delta)}.
$$
Thanks to \eqref{estimCas1sb} at $s=\frac12+\eta \delta$, this entails that (thanks once more to Proposition \ref{majBs21}) if in addition $\eta>0$ is so small that $\eta_0-\frac{\eta}2\geq k\eta_0$, then
\begin{multline}
 \|\de\|_{L_T^2 L^\infty} \leq \|\de\|_{L_T^2 \dot{B}_{2,1}^\frac32} \leq \|\de\|_{L_T^2 \dot{H}^{\frac32-\eta\delta}}^\frac12 \|\de\|_{L_T^2 \dot{H}^{\frac32+\eta\delta}}^\frac12\\
 \leq \DTt \ee^{\frac12 \left(\min(\aa_0, \frac23k'k''\eta_0 \delta) +\min(\aa_0, (\eta_0-\frac{\eta}2) \delta) \right)} \leq \DTt \ee^{\frac12 \left(\min(\aa_0, \frac23k'k''\eta_0 \delta) +\min(\aa_0, k\eta_0 \delta) \right)},
\end{multline}
and when we choose $k'=k''=\sqrt{k}<1$ we get that:
\begin{equation}
\|\de\|_{L_T^2 L^\infty} \leq \DTt \ee^{\frac12 \left(\min(\aa_0, \frac23k\eta_0 \delta) +\min(\aa_0, k\eta_0 \delta) \right)}.
\end{equation}
Finally, applying Proposition \ref{Estimdispnu} to $\Weh$ and also to $\Wei$ (as explained in the beginning of the present section) with $(d,p,r,q)=(0,2,\infty, 1)$ (and with the same arguments as in the previous section but with $\eta'=0$), we get that:
$$
\|\Weh+\Wei\|_{L_T^2 L^\infty} \leq \DTt \ee^{k\eta_0 \delta},
$$
gathering the last two estimates concludes the proof. $\blacksquare$

\subsection{On the optimality of the condition $\delta\leq \frac16$}

We wish to explain in this section why it is not possible (with the arguments of the present article) to improve the condition $\delta \leq \frac16$ into the one from \cite{FCcompl}: $\delta< \frac14$. According to the proof of Proposition \ref{estimStriPE} the estimates involving $L^6$ only require the condition $\delta\leq \frac13$, $\delta\leq \frac16$ being required by the ones involving $L^3$, and come from the estimates of $F_4, F_6, F_7$ and $F_{10}$. For some of them it is possible to improve the condition but the main limitation comes from $F_{10}$: we look for $k_1, k_2\in [2,\infty]$ such that $\frac1{k_1}+\frac1{k_2}=\frac12$ and,
$$
 |(F_{10}|\de)_{\dot{H}^s}| \leq \frac{\nu}{18}\|\de\|_{\dot{H}^{s+1}}^2 +\frac{C}{\nu^{\frac{s}{1-s}}} \|\Weh\|_{L^{k_1}}^{\frac2{1-s}} \|\de\|_{\dot{H}^s}^2 + C \|\n \Weh\|_{L^{k_2}}^2.
$$
At first sight it seems surprising that there is no better choice that $(k_1, k_2)=(6,3)$ and to understand this, let us focus on the Strichartz estimates involved by the previous bound. Choosing successively $(d,p,r,q)\in\{(0,\frac2{1-s},k_1,2), (1,2,k_2,2)\}$ we end-up with (for $\theta, \theta'\in[0,1]$) with the regularity indices:
$$
 \begin{cases}
 \sigma=\frac12+s-\frac3{k_1}+\frac{\theta}2 (1-\frac2{k_1}),\\
 \sigma'= \frac{3+\theta}{k_1}.
\end{cases}
$$
There exists $\theta, \theta'\in[0,1]$ such that $\sigma=\sigma'=\frac12+ \delta$ if, and only if, we have
\begin{equation}
 \frac1{k_1}\in \Big]\frac{s-\delta}3, \frac14 (s+\frac12-\delta)\Big] \cap \Big[\frac14(\frac12+\delta), \frac13(\frac12+\delta)\Big[.
 \label{intervk1}
\end{equation}
It is elementary to see that for any $a,b,c,d\in \R$ with $a<b$ and $c<d$,
$$
[a,b[ \cap ]c,d]\neq \emptyset \Longleftrightarrow a\leq d \mbox{ and } c<b,
$$
so the set in \eqref{intervk1} is nonempty if and only if $s\in [2\delta, \frac12+2\delta[$. In this configuration, the Strichartz estimates for $\Weh$ would write as follows:
\begin{equation}
\begin{cases}
 \|\Weh\|_{L_T^{\frac2{1-s}} L^{k_1}}\leq C \ee^{\frac12(\delta+\frac3{k_1}-s)} \|\Uoosc\|_{\dot{H}^{\frac12+\delta}},\\
 \|\n \Weh\|_{L_T^2 L^{k_2}}\leq C \ee^{\frac12(\frac12+\delta-\frac3{k_1})} \|\Uoosc\|_{\dot{H}^{\frac12+\delta}}.
\end{cases}
\label{estimStrik1}
\end{equation}
\begin{itemize}
 \item Finding $k_1$ when $s=\frac12$ and $\|\Uoosc\|_{\dot{H}^{\frac12+\delta}}\leq m(\ee)\ee^{-\frac{\delta}2}$ leads to:
$$
\begin{cases}
 \|\Weh\|_{L_T^4 L^{k_1}}\leq C \ee^{\frac12(\frac3{k_1}-\frac12)} m(\ee),\\
 \|\n \Weh\|_{L_T^2 L^{k_2}}\leq C \ee^{\frac12(\frac12-\frac3{k_1})} m(\ee),
\end{cases}
$$
and they are useful if both powers of $\ee$ are nonnegative, which leads to $k_1=6$.
\item Finding $k_1$ for any $s\in [\frac12, \frac12+\eta \delta]$ and $\|\Uoosc\|_{\dot{H}^{\frac12+\delta}}\leq \ee^{-\gamma}$ then requires that $[\frac12, \frac12+\eta \delta] \subset [2\delta, \frac12+2\delta[$ which is equivalent to $\delta\leq \frac14$ and $\eta<2$ (this one being true as $\eta \leq 2\eta_0 <1$). This leads to the estimates:
$$
\begin{cases}
 \|\Weh\|_{L_T^\frac2{1-s} L^{k_1}}\leq C \ee^{\frac12(\delta+\frac3{k_1}-s-2\gamma)},\\
 \|\n \Weh\|_{L_T^2 L^{k_2}}\leq C \ee^{\frac12(\frac12+\delta-\frac3{k_1}-2\gamma)},
\end{cases}
$$
Both exponents are positive for any $s\in [\frac12, \frac12+\eta \delta]$ if and only if:
$$
\frac1{k_1} \in \Big]\frac16 -\frac{2 \eta_0-\eta}3 \delta ,\frac16 +\frac23 \eta_0\delta \Big[,
$$
and we can put $\frac1{k_1}=\frac16 +\aa \delta$ with $\aa\in ]-\frac{2 \eta_0-\eta}3,\frac23 \eta_0[$ and $\theta, \theta'$ then write as follows
$$
(\theta, \theta')=\left( \frac{\frac12-s+\delta(1+3\aa)}{\frac13-\delta \aa}, \frac{\delta(1-3\aa)}{\frac16+\delta \aa}\right).
$$
As we already require that $\delta\leq \frac14$, both of them lie in $[0,1]$ if and only if $\aa \in[\frac14(1-\frac1{6\delta}), \frac14(\frac1{3\delta}-1)]$. Once more, the existence of such an $\aa$ is equivalent to the fact that:
$$
\Big[\frac14(1-\frac1{6\delta}), \frac14(\frac1{3\delta}-1)\Big] \cap \Big]-\frac{2 \eta_0-\eta}3,\frac23 \eta_0\Big[ \neq \emptyset,
$$
which is equivalent to
$$
\frac14(1-\frac1{6\delta}) < \frac23 \eta_0, \quad \mbox{and} \quad \frac14(\frac1{3\delta}-1)> -\frac{2 \eta_0-\eta}3.
$$
Both conditions are realized when $\delta\leq\frac16$. On the other hand if $\delta\in]\frac16, \frac14]$ the first condition means that $\eta_0>0$ is bounded from below by a positive constant and cannot be chosen as small as we wish. In other words, thanks to the definition of $\eta_0= \frac12(1-2\frac{\gamma}{\delta})$, the condition is equivalent to $\gamma<\frac18(\frac12 +\delta)$ which means $\gamma$ cannot be close to $\frac{\delta}2$ anymore (for example the condition becomes $\gamma<\frac3{32}$ when $\delta=\frac14$). So if we wish to choose $\gamma$ close to $\frac{\delta}2$ we need $\delta\leq \frac16$ and the only choice is $k_1=6$.
\end{itemize}

\section{Proof of Theorem \ref{ThRF}}

The proof will share the same steps as in the previous section, but keeping in mind that dealing with product of 2D and 3D functions will also induce a modification of the use of the Stichartz estimates (that will become anisotropic as in \cite{CDGG, CDGG2, CDGGbook}).

\subsection{Auxiliary systems}

Let us consider the initial data $\voe= \ub_0 + \woe$ with $\ub_0 \in [L^2(\R^2)]^3$ and $\woe \in [\dot{H}^\frac12(\R^3)\cap \dot{H}^{\frac12+\delta}(\R^3)]^3$ (both of them divergence-free). From the results recalled in the introduction:
\begin{itemize}
 \item there exists a global solution $\ub$ of System \eqref{NS2D},
 \item there exists a local strong solution $\we$ of System \eqref{PRF}, defined for some lifespan $\Te$ and for any $T<\Te$, $\we \in \dot{E}_T^\frac12$,
 \item moreover, the blow-up (or continuation) criterion is valid:
 $$
 \Te<\infty \Longrightarrow \int_0^{\Te} \|\n \we(t)\|_{\dot{H}^\frac12}^2 dt =\infty,
 $$
 \item finally, as $\woe\in \dot{H}^\frac12\cap \dot{H}^{\frac12+\delta}$ then for all $T<\Te$, and $s\in[\frac12, \frac12+\delta]$, $\we \in \dot{E}_T^s$. 
\end{itemize}
Introducing $\We$ as the global solution of the following linear system:
\begin{equation}
 \begin{cases}
  \d_t \We -\nu \D \We +\frac{1}{\ee}\mathbb{P} (e_3\wedge \We)=0,\\
  {\We}_{|t=0}=\woe,
 \end{cases}
\label{LRFb}\tag{$LRF_\ee$}
\end{equation}
we define on $\de\overset{def}{=} \ve-\ub-\We =\we-\We$, which satisfies:
\begin{equation}
 \begin{cases}
  \d_t \de -\nu \D \de +\frac{1}{\ee} \mathbb{P} (e_3 \wedge \de) = \Sum_{i=1}^{8} G_i,\\
  {\de}_{|t=0}= 0,
 \end{cases}
 \label{SystdeRF}
\end{equation}
with:
\begin{equation}
\begin{cases}
 G_1 \overset{def}{=}-\mathbb{P}(\de \cdot \n \de), \qquad G_2 \overset{def}{=}-\mathbb{P}(\de \cdot \n \We), \qquad G_3 \overset{def}{=}-\mathbb{P}(\We \cdot \n \de),\\
 G_4 \overset{def}{=}-\mathbb{P}(\We \cdot \n \We),\qquad G_5 \overset{def}{=}-\mathbb{P}(\de \cdot \n \ub),\qquad G_6 \overset{def}{=}-\mathbb{P}(\ub \cdot \n \de),\\
 G_7 \overset{def}{=}-\mathbb{P}(\We \cdot \n \ub),\qquad G_8 \overset{def}{=}-\mathbb{P}(\ub \cdot \n \We).
\end{cases}
 \label{G18}
\end{equation}

\subsection{Estimates on $\de$}
Let us assume that $\ee$ satisfies \eqref{condeps1} and \eqref{condeps3} and assume by contradiction that $\Te<\infty$, then by the continuation criterion, we have:
\begin{equation}
 \int_0^{\Te} \|\n \we(t)\|_{\dot{H}^\frac12}^2 dt = \infty,
 \label{critereexpl2}
\end{equation}
If we put (where the constant $C$ refers to the one from \eqref{EstimG1234})
\begin{equation}
 T_\ee \overset{def}{=} \sup \{t\in[0,\Te[, \quad \forall t'\leq t, \|\de(t')\|_{\dot{H}^\frac12} \leq \frac{\nu}{4C}\},
\end{equation}
As $\delta_\ee(0)=0$ then $T_\ee>0$. Now assume by contradiction that:
\begin{equation}
T_\ee<\Te
 \label{Contrad3},
\end{equation}
and the $\dot{H}^s$ innerproduct of \eqref{SystdeRF} with $\de$ leads to:
$$
\frac12 \frac{d}{dt} \|\de(t)\|_{\dot{H}^s}^2 +\nu \|\n \de(t)\|_{\dot{H}^s}^2 \leq \Sum_{j=1}^8 (G_j|\de)_{\dot{H}^s}.
$$
As the method is similar to what we did previously, we will skip details about the following terms whose estimates are done as in the first section:
\begin{equation}
 \begin{cases}
 \vspace{0.1cm}
  |(G_1|\de)_{\dot{H}^s}| \leq C \|\de\|_{\dot{H}^\frac12} \|\de\|_{\dot{H}^{s+1}}^2,\\
  \vspace{0.1cm}
  |(G_2|\de)_{\dot{H}^s}| \leq \frac{\nu}{14}\|\de\|_{\dot{H}^{s+1}}^2 +\frac{C}{\nu} \|\n \We\|_{L^3}^2\|\de\|_{\dot{H}^s}^2,\\
  |(G_3|\de)_{\dot{H}^s}| \leq \frac{\nu}{14}\|\de\|_{\dot{H}^{s+1}}^2 +\frac{C}{\nu^3} \|\We\|_{L^6}^4 \|\de\|_{\dot{H}^s}^2,\\
  |(G_4|\de)_{\dot{H}^s}| \leq \frac{\nu}{14}\|\de\|_{\dot{H}^{s+1}}^2 +\frac{C}{\nu^{\frac{s}{1-s}}} \|\We\|_{L^6}^{\frac2{1-s}} \|\de\|_{\dot{H}^s}^2 + C \|\n \We\|_{L^3}^2.
 \end{cases}
\label{EstimG1234}
\end{equation}
We will only focus on what changes, namely the terms involving products of 2D and 3D functions. The first two terms are easily estimated with the usual arguments thanks to Proposition \ref{prod2D3D}, as $s\in[\frac12, \frac12+\eta \delta]\subset[0,1[$:
\begin{multline}
 |(G_5|\de)_{\dot{H}^s}| \leq \|\de \cdot \n \ub\|_{\dot{H}^{s-1}} \|\de\|_{\dot{H}^{s+1}} \leq \|\de\|_{\dot{H}^s} \|\n \ub\|_{\dot{H}^0} \|\de\|_{\dot{H}^{s+1}}\\
\leq \frac{\nu}{14}\|\de\|_{\dot{H}^{s+1}}^2 +\frac{C}{\nu} \|\ub\|_{\dot{H}^1}^2 \|\de\|_{\dot{H}^s}^2
\label{EstimG5}
\end{multline}
Similarly, we easily get (with $(s_1,s_2)=(\frac12, s-\frac12)$):
\begin{equation}
|(G_6|\de)_{\dot{H}^s}| \leq \frac{\nu}{14}\|\de\|_{\dot{H}^{s+1}}^2 +\frac{C}{\nu^3} \|\ub\|_{L^2}^2 \|\ub\|_{\dot{H}^1}^2 \|\de\|_{\dot{H}^s}^2.
 \label{EstimG6}
\end{equation}
Now we can turn to the last terms and obtain, adapting the arguments from the previous section ($s\in[0,1]$), that: 
\begin{multline}
 |(G_7|\de)_{\dot{H}^s}| \leq \|\We \cdot \n \ub\|_{L^2} \|\de\|_{\dot{H}^{2s}} \leq C \|\We\|_{L_{h,v}^{\infty, 2}} \|\n \ub\|_{L^2(\R^2)} \|\de\|_{\dot{H}^s}^{1-s} \|\de\|_{\dot{H}^{s+1}}^s\\
\leq \frac{\nu}{14}\|\de\|_{\dot{H}^{s+1}}^2 +\frac{C}{\nu^\frac{s}{1-s}} \|\ub\|_{\dot{H}^1}^2 \|\de\|_{\dot{H}^s}^2 +\|\ub\|_{\dot{H}^1}^{2s} \|\We\|_{L_{h,v}^{\infty, 2}}^2,
\label{EstimG7}
\end{multline}
and thanks to the Sobolev injection $\dot{H}^\frac12(\R^2) \hookrightarrow L^4(\R^2)$ , and the Young inequality with $(2,\frac2{1-s},\frac2{s})$:
\begin{multline}
 |(G_8|\de)_{\dot{H}^s}| \leq \|\ub \cdot \n \We\|_{L^2} \|\de\|_{\dot{H}^{2s}} \leq C\|\ub\|_{L^4(\R^2)} \|\n \We\|_{L_{h,v}^{4, 2}} \|\de\|_{\dot{H}^{2s}}\\
 \leq C \left(\|\ub\|_{L^2(\R^2)}^\frac12 \|\ub\|_{\dot{H}^1(\R^2)}^{s-\frac12} \|\n \We\|_{L_{h,v}^{4, 2}}\right) \left(\|\ub\|_{\dot{H}^1(\R^2)} \|\de\|_{\dot{H}^s}\right)^{1-s} \|\de\|_{\dot{H}^{s+1}}^s\\
\leq \frac{\nu}{14}\|\de\|_{\dot{H}^{s+1}}^2 +\frac{C}{\nu^\frac{s}{1-s}} \|\ub\|_{\dot{H}^1}^2 \|\de\|_{\dot{H}^s}^2 +\|\ub\|_{L^2} \|\ub\|_{\dot{H}^1}^{2s-1} \|\n \We\|_{L_{h,v}^{4, 2}}^2.
\label{EstimG8}
\end{multline}
As in the previous section, collecting \eqref{EstimG1234} to \eqref{EstimG8}, and thanks to the energy equality from Theorem \ref{thN2D}, there exists some constant $\Bo=\Bo(\nu, s, \|\ub_0\|_{L^2})>0$ such that we can write that for any $t\leq T_\ee$,
\begin{multline}
 \|\de(t)\|_{\dot{H}^s}^2 +\frac{\nu}{2} \int_0^t \|\n \de(\tau)\|_{\dot{H}^s}^2 d\tau \leq \Bo e^{\Bo \left(1 +\|\n \We\|_{L^2 L^3}^2 +\|\We\|_{L^4 L^6}^4 +\|\We\|_{L^{\frac2{1-s}} L^6}^{\frac2{1-s}}\right)}\\
 \times \left(\|\n \We\|_{L^2 L^3}^2 +\|\We\|_{L^{\frac2{1-s}}L_{h,v}^{\infty, 2}}^2 +\|\n \We\|_{L^{\frac2{\frac32-s}} L_{h,v}^{4, 2}}^2 \right).
\label{estimCas1RF}
\end{multline}

\subsection{Proof of Point 1}

We can now plug in \eqref{estimCas1RF} the Strichartz estimates from Proposition \ref{estimStriRF} and obtain that for any $k\in]0,1[$ (as close to 1 as wished) fixed, any $s\in[\frac12, \frac12+\eta \delta]$, there exists a constant $\Bo$ such that for all $t\leq T_\ee$,
\begin{multline}
 \|\de(t)\|_{\dot{H}^s}^2 +\frac{\nu}{2} \int_0^t \|\n \de(\tau)\|_{\dot{H}^s}^2 d\tau\\
 \leq \Bo e^{\Bo \left(1 +\ee^{2\eta_0 \delta} +\ee^{4\eta_0 \delta} +\ee^{\frac1{1-s}(\frac12+2\eta_0\delta-s)}\right)} \times \left(\ee^{2\eta_0 \delta} + \ee^{k(\frac12+2\eta_0 \delta-s)}\right).
\end{multline}
Now if $\ee>0$ is so small that \eqref{condeps1} is true, then putting $\Do=\Bo e^{2\Bo}$:
\begin{equation}
 \|\de(t)\|_{\dot{H}^s}^2 +\frac{\nu}{2} \int_0^t \|\n \de(\tau)\|_{\dot{H}^s}^2 d\tau \leq \Do \ee^{k(\frac12+2\eta_0 \delta-s)}.
 \label{estimdeltaRF2}
\end{equation}
From this the rest of the boostrap argument is classic and similar to what is done in \cite{FCPAA, FCcompl}: assuming that $\ee$ is so small that (taking $s=\frac12$):
\begin{equation}
 \Do \ee^{2k\eta_0 \delta} \leq \left(\frac{\nu}{8C}\right)^2,
 \label{condeps3}
\end{equation}
then we obtain that for all $t\leq T_\ee$, $\|\de(t)\|_{\dot{H}^\frac12} \leq \frac{\nu}{8C}$, which contradicts the definition of $T_\ee$ and the fact that $\|\de(T_\ee)\|_{\dot{H}^\frac12}=\frac{\nu}{4C}$. Then $T_\ee=\Te$, so that, as $\we=\de+\We$, we obtain that for all $t<\Te$,
$$
\int_0^t \|\we(\tau)\|_{\dot{H}^\frac12}^2 d\tau \leq \int_0^t \|\de(\tau)\|_{\dot{H}^\frac12}^2 d\tau +\int_0^t \|\We(\tau)\|_{\dot{H}^\frac12}^2 d\tau \leq \frac1{\nu}\left(\Do \ee^{2k\eta_0 \delta} +\|\woe\|_{\dot{H}^\frac12}^2\right).
$$
Even if we only control the norm of $\woe$ in $\dot{H}^{\frac12+ c\delta}\cap \dot{H}^{\frac12 + \delta}$ the previous quantity is finite, which contradicts \eqref{critereexpl2}, so that $\Te=\infty$.

\subsection{Proof of Point 2}

This part is nearly identical to what we did in Section \ref{Point2b} so we will not give much details: for any $\eta\in[0,2\eta_0[$, any $k\in]0,1[$ (near to $1$) and any $\eta'\in[0,\eta[$, from \eqref{estimdeltaRF2} with $s\in\{\frac12, \frac12+\eta\delta\}$ we get:
\begin{multline}
 \||D|^{\eta'\delta} \de\|_{L^2 L^\infty} \leq \||D|^{\eta'\delta} \de\|_{L^2 \dot{B}_{2,1}^\frac32} \leq \||D|^{\eta'\delta} \de\|_{L^2 \dot{H}^{\frac32-\eta'\delta}}^{1-\frac{\eta'}{\eta}} \||D|^{\eta'\delta} \de\|_{L^2 \dot{H}^{\frac32+(\eta-\eta')\delta}}^{\frac{\eta'}{\eta}}\\
 \leq \|\de\|_{L^2 \dot{H}^\frac32}^{1-\frac{\eta'}{\eta}} \|\de\|_{L^2 \dot{H}^{\frac32+\eta\delta}}^{\frac{\eta'}{\eta}} \leq \Do \ee^{k\delta \left(\eta_0(1-\frac{\eta'}{\eta})  +(\eta_0-\frac12\eta)\frac{\eta'}{\eta}\right)} =\Do \ee^{k\delta (\eta_0-\frac12\eta')}
\end{multline}
Thanks to Proposition \ref{estimStriRF} with $(d,p,m,q)=(\eta'\delta, 2, \infty, 1)$ and doing the same as in Section \ref{Point2b} we obtain:
$$
\||D|^{\eta'\delta} \We\|_{L^2 L^\infty} \leq \Co C_{k,\eta_0, \eta', \nu} \ee^{k\delta (\eta_0-\frac12\eta')},
$$
which ends the proof. $\blacksquare$

\subsection{Proof of Point 3}

Let us fix some $k,k',k''\in]0,1[$. First thanks to \eqref{estimdeltaRF2} for $s=\frac12+\eta\delta$, we get that:
\begin{equation}
 \|\de\|_{L^2 \dot{H}^{\frac32+\eta \delta}} \leq \Do \ee^{k'(\eta_0-\frac{\eta}2) \delta}.
 \label{de1}
\end{equation}
For the same reason as in Section \ref{Point2c}, the estimates for $G_8$ has to be changed when $s<\frac12$ in a similar way as we did for $F_6$ and $F_7$, and thanks to Proposition \ref{prod2D3D}, we get:
\begin{multline}
 |(G_8|\de)_{\dot{H}^s}| \leq \|\ub \cdot \n \We\|_{L^2} \|\de\|_{\dot{H}^{2s}} \leq C\|\ub\|_{L^2(\R^2)} \|\n \We\|_{L_{h,v}^{\infty, 2}} \|\de\|_{\dot{H}^{2s}}\\
 \leq C \left(\|\ub\|_{L^2(\R^2)} \|\n \We\|_{L_{h,v}^{\infty, 2}}^\frac1{2-s}\right) \left( \|\n \We\|_{L_{h,v}^{\infty, 2}}^\frac{1-s}{2-s}\|\de\|_{\dot{H}^s}^{1-s}\right) \|\de\|_{\dot{H}^{s+1}}^s\\
\leq \frac{\nu}{14}\|\de\|_{\dot{H}^{s+1}}^2 +\frac{C}{\nu^\frac{s}{1-s}} \|\n \We\|_{L_{h,v}^{\infty, 2}}^\frac2{2-s} \|\de\|_{\dot{H}^s}^2 +\|\ub\|_{L^2}^2 \|\n \We\|_{L_{h,v}^{\infty, 2}}^\frac2{2-s}.
\label{EstimG8bis}
\end{multline}
So that \eqref{estimCas1RF} turns into
\begin{multline}
 \|\de(t)\|_{\dot{H}^s}^2 +\frac{\nu}{2} \int_0^t \|\n \de(\tau)\|_{\dot{H}^s}^2 d\tau \leq \Bo e^{\Bo \left(1 +\|\n \We\|_{L^2 L^3}^2 +\|\We\|_{L^4 L^6}^4 +\|\We\|_{L^{\frac2{1-s}} L^6}^{\frac2{1-s}} +\|\n \We\|_{L^\frac2{2-s} L_{h,v}^{\infty, 2}}^\frac2{2-s}\right)}\\
 \times \left(\|\n \We\|_{L^2 L^3}^2 +\|\We\|_{L^{\frac2{1-s}}L_{h,v}^{\infty, 2}}^2 +\|\n \We\|_{L^\frac2{2-s} L_{h,v}^{\infty, 2}}^\frac2{2-s} \right).
\end{multline}
Now with the same arguments we obtain that:
\begin{equation}
\|\de\|_{L^2 \dot{H}^{\frac32-\eta \delta}}^2 \leq \Do \left(\ee^{2\eta_0 \delta} +\ee^{k'(2\eta_0+\eta) \delta} +\ee^{k'\frac{2\eta_0+\eta}{\frac32+\eta\delta} \delta} \right).
\label{de2}
\end{equation}
If we choose $\eta>0$ so small that
$$
\eta_0-\frac{\eta}2 \geq k'' \eta_0,\quad \mbox{and} \quad \frac{2\eta_0+\eta}{\frac32+\eta\delta} \geq \frac43 k''\eta_0,
$$
then thanks to Proposition \ref{estimStriRF}, we obtain (with $k'=k''=\sqrt{k}$):
$$
 \|\de\|_{L^2 L^\infty} \leq \|\de\|_{L^2 \dot{B}_{2,1}^\frac32} \leq \|\de\|_{L^2 \dot{H}^{\frac32-\eta\delta}}^\frac12 \|\de\|_{L^2 \dot{H}^{\frac32+\eta\delta}}^\frac12\leq \Do \ee^{\frac56 k\eta_0\delta}.
$$
Then, using Proposition \ref{estimStriRF} with $(d,p,m,q)=(0,2,\infty,1)$, we obtain:
$$
\|\We\|_{L^2 L^\infty} \leq \Do \ee^{k\eta_0 \delta},
$$
which concludes the proof. $\blacksquare$

\section{Appendix}

\subsection{Notations, Sobolev spaces and Littlewood-Paley decomposition}

We refer to the appendix of \cite{FCPAA} for general notations and properties of the Sobolev spaces and the Littlewood-Paley decomposition (together with the classical properties). For a complete presentation, we refer to \cite{Dbook}. Let us first mention the following lemma whose proof is close to Lemma $5$ from \cite{FCestimLp}  (see also Section 2.11 in \cite{Dbook}):
\begin{lem} \label{majBs21}
 \sl{For any $\aa, \beta>0$ there exists a constant $C_{\aa, \beta}>0$ such that for any $u\in \dot{H}^{s-\aa} \cap \dot{H}^{s+\beta}$, then $u\in\dot{B}_{2,1}^s$ and:
\begin{equation}
 \|u\|_{\dot{B}_{2,1}^s} \leq C_{\aa, \beta} \|u\|_{\dot{H}^{s-\aa}}^{\frac{\beta}{\aa + \beta}} \|u\|_{\dot{H}^{s+\beta}}^{\frac{\aa}{\aa + \beta}}.
\end{equation}
 }
\end{lem}
\begin{prop}
 \sl{\cite{Dbook} We have the following continuous injections:
$$
 \begin{cases}
\mbox{For any } p\geq 1, & \dot{B}_{p,1}^0 \hookrightarrow L^p,\\
\mbox{For any } p\in[2,\infty[, & \dot{B}_{p,2}^0 \hookrightarrow L^p,\\
\mbox{For any } p\in[1,2], & \dot{B}_{p,p}^0 \hookrightarrow L^p.
\end{cases}
$$
}
 \label{injectionLr}
\end{prop}
An alternative to the classical $L_t^p \dot{B}_{q,r}^s$-type estimates is provided by the Chemin-Lerner time-space Besov spaces: as explained in the following definition, the integration in time is performed before the summation with respect to the frequency decomposition index:
\begin{defi} \cite{Dbook}
 \sl{For $s,t\in \R$ and $a,b,c\in[1,\infty]$, we define the following norm
 $$
 \|u\|_{\tilde{L}_t^a \dot{B}_{b,c}^s}= \Big\| \left(2^{js}\|\ddj u\|_{L_t^a L^b}\right)_{j\in \Z}\Big\|_{l^c(\Z)}.
 $$
 The space $\tilde{L}_t^a \dot{B}_{b,c}^s$ is defined as the set of tempered distributions $u$ such that $\lim_{j \rightarrow -\infty} S_j u=0$ in $L^a([0,t],L^\infty(\R^d))$ and $\|u\|_{\tilde{L}_t^a \dot{B}_{b,c}^s} <\infty$.
 }
 \label{deftilde}
\end{defi}
We refer once more to \cite{Dbook} (Section 2.6.3) for more details and will only recall the following proposition:
\begin{prop}
\sl{
For all $a,b,c\in [1,\infty]$ and $s\in \R$:
     $$
     \begin{cases}
    \mbox{if } a\leq c,& \forall u\in L_t^a \dot{B}_{b,c}^s, \quad \|u\|_{\tilde{L}_t^a \dot{B}_{b,c}^s} \leq \|u\|_{L_t^a \dot{B}_{b,c}^s}\\
    \mbox{if } a\geq c,& \forall u\in\tilde{L}_t^a \dot{B}_{b,c}^s, \quad \|u\|_{\tilde{L}_t^a \dot{B}_{b,c}^s} \geq \|u\|_{L_t^a \dot{B}_{b,c}^s}.
     \end{cases}
     $$
     \label{Propermut}
     }
\end{prop}

\subsection{Strichartz estimates for the primitive system}
\label{Disp}

\subsubsection{Statements of the results}

Consider the following system (in the case $\nu=\nu'$, we have $L=\nu\Delta$):
\begin{equation}
\begin{cases}
\d_t f-(\nu\Delta-\frac{1}{\ee} \mathbb{P} \cA) f=\Fe,\\
f_{|t=0}=f_0.
\end{cases}
\label{systdisp}
\end{equation}
Let us recall the Strichartz estimates obtained in \cite{FCcompl} (we refer to \cite{FCPAA, FCcompl} for details about the system, its analysis as well as the notations used).

\begin{prop}
 \sl{For any $d\in \R$, $r\geq 2$, $q\geq 1$, $\theta\in[0,1]$ and $p\in[1, \frac{4}{\theta (1-\frac{2}{r})}]$, there exists a constant $C=C_{F, p,\theta,r}$ such that for any $f$ solving \eqref{systdisp} for initial data $f_0$ and external force $\Fe$ both with zero divergence and potential vorticity, then
 \begin{equation}
  \||D|^d f\|_{\tilde{L}_t^p\dot{B}_{r, q}^0} \leq \frac{C_{F, p,\theta,r}}{\nu^{\frac{1}{p}-\frac{\theta}{4}(1-\frac{2}{r})}} \ee^{\frac{\theta}{4}(1-\frac{2}{r})} \left( \|f_0\|_{\dot{B}_{2, q}^{\sigma_1}} +\|\Fe\|_{\tilde{L}_t^1 \dot{B}_{2, q}^{\sigma_1}} \right),
 \end{equation}
 where $\sigma_1= d+\frac32-\frac{3}{r}-\frac{2}{p}+\frac{\theta}{2} (1-\frac{2}{r})$.
\label{Estimdispnu}
 }
\end{prop}
Let us first state the following modified Strichartz estimates needed to fit to the regularity of the external force term G (see \eqref{G}) under the actual assumptions on $\tUqg$.
\begin{prop}
 \sl{For any $d\in \R$, $k\in ]1,2]$, $r\geq 2$, $q\geq 1$, $\theta\in[0,1]$ and $p\in[2, \frac{4}{\theta (1-\frac{2}{r})}[$, there exists a constant $C=C_{F, p,\theta,r,k}$ such that for any $f$ solving \eqref{systdisp} with zero initial data and an external force $\Fe$  with zero divergence and potential vorticity, then
 \begin{equation}
  \||D|^d f\|_{\tilde{L}_t^p\dot{B}_{r, q}^0} \leq \frac{C_{F, p,\theta,r,k}}{\nu^{1-\frac1{k}+\frac{1}{p}-\frac{\theta}{4}(1-\frac{2}{r})}} \ee^{\frac{\theta}{4}(1-\frac{2}{r})} \|\Fe\|_{\tilde{L}_t^k \dot{B}_{2, q}^{\sigma_2}},
 \end{equation}
 where $\sigma_2= d-\frac12+\frac2{k}-\frac{3}{r}-\frac{2}{p}+\frac{\theta}{2} (1-\frac{2}{r})$.
\label{Estimdispnu2}
 }
\end{prop}

\begin{rem}
\sl{The case $k=1$ is not covered by the second result but is dealt with in the first one. The condition on the $p$-index is more restrictive.}
\end{rem}

Now, as a consequence of Propositions \ref{Estimdispnu}, \ref{Estimdispnu2} and \ref{PropestimG}, we can bound the various terms from \eqref{estimCas1} involving $\Weh$ and $\Wei$ . We collect these estimates in the following proposition:

\begin{prop}
 \sl{Under the previous notations, if $\delta\leq \frac16$, for any $s\in[\frac12, \frac12+\eta \delta]$, there exists a constant $C=C(F,\delta,s)>0$ such that:
 $$
 \begin{cases}
  \displaystyle{\|\Weh\|_{L_T^4 L^6} +\nu^{\frac14} \|\n \Weh\|_{L_T^2 L^3} \leq \frac{C}{\nu^{\frac{1-2\delta}4}} \ee^{\frac{\delta}2} \|\Uoosc\|_{\dot{H}^{\frac12+\delta}}}\\
  \displaystyle{\|\Weh\|_{L_T^{\frac2{1-s}}L^6} +\nu^{\frac14} \|\n \Weh\|_{L_T^{\frac2{\frac32-s}}L^3} \leq\frac{C}{\nu^{\frac{1-2\delta}4}} \ee^{\frac12(\frac12+\delta-s)} \|\Uoosc\|_{\dot{H}^{\frac12+\delta}}}\\
  \displaystyle{\|\Wei\|_{L_T^{\frac2{1-s}}L^6} +\nu^{\frac14} \|\n \Wei\|_{L_T^{\frac2{\frac32-s}}L^3} \leq\frac{C}{\nu^{\frac{5-2\delta}4}} \ee^{\frac12(\frac12+\delta-s)} \Co^2 e^{\frac{C}{\nu}\CT}.}
 \end{cases}
$$
Under Assumption $(H_5)$, if $\delta\leq \frac16$, for any $k\in]0,1[$ and $s\in[\frac12-\eta \delta, \frac12]$ with $\eta\leq 2\eta_0$, there exists a constant $C=C_{F,\delta,s,k,c,\gamma}>0$ such that the previous estimates remain true, except those involving the $L_T^{\frac2{\frac32-s}}L^3$-norms, which are replaced by:
$$
\|\n \Weh\|_{L_T^\frac2{2-s}L^\frac3{s}} +\nu \|\n \Wei\|_{L_T^\frac2{2-s}L^\frac3{s}}\leq \frac{C}{\nu^{1-\frac{s}2-\gamma-k\delta\eta_0}} \ee^{k\eta_0 \delta}(1+ \ee^{\gamma} e^{\frac{C}\nu \CT}).
$$
 }
 \label{estimStriPE}
\end{prop}

\subsubsection{Proof of Proposition \ref{Estimdispnu2}}
\label{preuvedispPE}
As we outlined in \cite{FCPAA, FCcompl}, for any divergence-free and with zero potential vorticity initial data $g_0$, the operators $\cP$ and $\cP_{3+4}$ (on one hand), $\cQ$ and $\cP_{2}$ (on the other hand) coincide when $\nu=\nu'$ (see the cited articles for precisions and notations): 
$$
g_0=\mathbb{P} g_0=\cP \mathbb{P} g_0= \mathbb{P}_{3+4} \mathbb{P} g_0 =\mathbb{P}_{3+4} g_0.
$$
We will denote as $S_\ee(t)$ the associated semi-group, that is $g(t)=S_\ee(t)g_0$ is the unique solution of \eqref{systdisp} with initial data $g_0$ and no external force: 
$$
S_\ee(t)g_0 =g(t)=\mathcal{F}^{-1}\left(e^{-\nu t |\xi|^2 +i\frac{t}{\ee} \frac{|\xi|_F}{F|\xi|}} \mathcal{P}_3(\xi,\ee) \hat{g_0}(\xi) +e^{-\nu t |\xi|^2 -i\frac{t}{\ee} \frac{|\xi|_F}{F|\xi|}} \mathcal{P}_4(\xi,\ee) \hat{g_0}(\xi) \right),
$$
and in order to simplify we will write:
$$
S_\ee(t)g_0=\mathcal{F}^{-1}\left(e^{-\nu t |\xi|^2 +i\frac{t}{\ee} \frac{|\xi|_F}{F|\xi|}} \hat{g_0}(\xi) \right),
$$
So that, thanks to the Duhamel formula, the solution $f$ from Proposition \ref{Estimdispnu2} writes:
$$
f(t,x)= \int_0^t S_\ee(t-\tau)\Fe (\tau,x)d\tau.
$$
We will only focus on what is new (and refer to \cite{FCPAA, FCcompl} for details or notations). If $\varphi$ is the usual truncation function involved in the Littlewood-Paley decomposition, let us denote by $\varphi_1$ another smooth truncation function, with support in a slightly larger annulus than $\mbox{supp }\varphi$ (say for example the annulus centered at zero and of radii $\frac12$ and $3$) and equal to $1$ on $\mbox{supp }\varphi$. For given $p,r \geq 1$, let $\cB$ be the set:
$$
\cB\overset{def}{=}\{\psi \in \cC_0^\infty (\R_+\times \R^3, \R), \quad \|\psi\|_{L^{\bar{p}}(\R_+, L^{\bar{r}}(\R^3))}\leq 1\},
$$
then we follow the same classical steps, for any $j\in \Z$:
\begin{multline}
 \|\ddj f\|_{L^p L^r}= \sup_{\psi \in \cB} \int_0^\infty \int_{\R^3} \left(\int_0^t S_\ee(t-\tau)\ddj\Fe (\tau,x)d\tau\right) \psi(t,x) dx dt\\
 =C \sup_{\psi \in \cB} \int_0^\infty \int_{\R^3} \left(\int_0^\infty e^{-\nu (t-\tau)|\xi|^2+i\frac{t-\tau}{\ee}\frac{|\xi|_F}{F|\xi|}} \varphi_1(2^{-j} \xi) \hat{\psi}(t,\xi) \textbf{1}_{\{\tau \leq t\}} dt\right) \widehat{\ddj \Fe}(\tau, \xi) d\xi d\tau\\
 \leq C \sup_{\psi \in \cB} \int_0^\infty \|\widehat{\ddj \Fe}(\tau,.)\|_{L^2}\\
 \times \left(\int_{\R^3} \int_0^\infty \int_0^\infty e^{-\nu (t+t'-2\tau)|\xi|^2+i\frac{t-t'}{\ee}\frac{|\xi|_F}{F|\xi|}} \varphi_1(2^{-j} \xi)^2 \hat{\psi}(t,\xi) \overline{\hat{\psi}(t',\xi)} \textbf{1}_{\{\tau \leq t\}} \textbf{1}_{\{\tau \leq t'\}} dt dt'd\xi\right)^{\frac12} d\tau\\
 \leq C \sup_{\psi \in \cB} \int_0^\infty \|\ddj \Fe(\tau,.)\|_{L^2}\\
 \times \left(\int_0^\infty \int_0^\infty \|L_j(\frac{t-t'}{\ee})\psi(t,.)\|_{L^r} \|e^{\nu(t+t'-2\tau)\Delta} \varphi_1(2^{-j}D) \overline{\psi(t',.)}\|_{L^{\bar{r}}} \textbf{1}_{\{\tau \leq \min(t,t')\}} dtdt'\right)^{\frac12} d\tau,
\label{estimTT1}
 \end{multline}
with $L_j(\sigma)$ defined as in \cite{FCcompl}:
$$
 L_j(\sigma)g=\int_{\R^3} e^{ix\cdot \xi +i \sigma \frac{|\xi|_F}{F|\xi|}} \varphi_1(2^{-j}|\xi|) \hat{g}(\xi) d\xi.
$$
We refer to \cite{FCcompl}, for the proof that for all $r\in[2,\infty]$ and $\theta \in[0,1]$:
$$
\begin{cases}
 \|e^{\sigma\Delta} \varphi_1(2^{-j}D) g\|_{L^{\bar{r}}} \leq C'e^{-\frac{\sigma}4 2^{2j}} \|g\|_{L^{\bar{r}}}\\
 \|L_j(\sigma)g\|_{L^r} \leq (C_F)^{1-\frac2{r}} \frac{2^{3j(1-\frac2{r})}}{|\sigma|^{\frac{\theta}2 (1-\frac2{r})}} \|g\|_{L^{\bar{r}}}.
\end{cases}
$$
Going back to \eqref{estimTT1}, we get:
\begin{multline}
 \|\ddj f\|_{L^p L^r} \leq C_F^{\frac12-\frac1{r}} 2^{3j(\frac12-\frac1{r})} \ee^{\frac{\theta}4 (1-\frac2{r})} \sup_{\psi \in \cB}\int_0^\infty \|\ddj \Fe(\tau,.)\|_{L^2} \mathbb{K}(\tau)^\frac12 d\tau\\
\leq C_F^{\frac12-\frac1{r}} 2^{3j(\frac12-\frac1{r})} \ee^{\frac{\theta}4 (1-\frac2{r})} \|\ddj \Fe\|_{L^k L^2} \times \sup_{\psi \in \cB} \left(\int_0^\infty  \mathbb{K}(\tau)^{\frac{\bar{k}}2} d\tau \right)^{\displaystyle{\frac1{\bar{k}}}},
\label{estimJen1}
\end{multline}
where
$$
\mathbb{K}(\tau) \overset{def}{=} \int_0^\infty \int_0^\infty e^{-\frac{\nu}4 2^{2j} (t+t'-2\tau)} \frac{\|\psi(t')\|_{L^{\bar{r}}} \|\psi(t)\|_{L^{\bar{r}}}}{|t-t'|^{\frac{\theta}2 (1-\frac2{r})}} \textbf{1}_{\{\tau \leq \min(t,t')\}} dtdt'.
$$
Next we use Jensen's inequality in the following formulation (we refer to \cite{BoyerFabrie}, Prop. $II.2.20$):
\begin{prop}
 \sl{Let $\Omega\subset \R^d$ be an open set, and $\eta \in L^1(\Omega)$ a nonnegative function. For any function $f$ such that $|f|^\aa \eta \in L^1(\Omega)$ for some $\aa \in [1,\infty[$, we have $f\eta \in L^1(\Omega)$ and
 $$
 \Big|\int_{\Omega}f\eta dx\Big|^\aa \leq \|\eta\|_{L^1}^{\aa-1} \int_{\Omega} |f|^\aa \eta dx.
 $$}
\end{prop}
Choosing $\aa=\frac{\bar{k}}2$, $\Omega=]0, \infty[^2$, $f(t,t')=e^{-\frac{\nu}8 2^{2j} (t+t'-2\tau)} \textbf{1}_{\{\tau \leq \min(t,t')\}}$ and,
$$
\eta(t,t')= \frac{f_\tau(t) f_\tau(t')}{|t-t'|^{\frac{\theta}2 (1-\frac2{r})}}, \quad \mbox{with }f_\tau (t)=e^{-\frac{\nu}8 2^{2j} (t-\tau)} \textbf{1}_{\{\tau \leq t\}} \|\psi(t)\|_{L^{\bar{r}}},
$$
we obtain that:
\begin{multline}
  \mathbb{K}(\tau)^{\frac{\bar{k}}2} \leq \left(\int_0^\infty \int_0^\infty \frac{f_\tau(t) f_\tau(t')}{|t-t'|^{\frac{\theta}2 (1-\frac2{r})}} dt dt'\right)^{\frac{\bar{k}}2-1}\\
  \times \int_0^\infty \int_0^\infty e^{-\frac{\nu}8 (1+\frac{\bar{k}}2) 2^{2j} (t+t'-2\tau)} \frac{\|\psi(t')\|_{L^{\bar{r}}} \|\psi(t)\|_{L^{\bar{r}}}}{|t-t'|^{\frac{\theta}2 (1-\frac2{r})}} \textbf{1}_{\{\tau \leq \min(t,t')\}} dt dt',
\label{estimJen2}
  \end{multline}
\begin{rem}
 \sl{This is here that we require $\frac{\bar{k}}2 \in [1,\infty[$, that is $k\in ]1,2]$.}
\end{rem}
The first integral is dealt with the Hardy-Littlewood-Sobolev estimates as in \cite{FCPAA, FCcompl}: introducing $\frac1{q_1}=1-\frac{\theta}4 (1-\frac2{r})$ (which is in $[1, \infty[$), and some constant $\mathbb{E}$ (depending on $\theta,r$)
\begin{multline}
\int_0^\infty \int_0^\infty \frac{f_\tau(t) f_\tau(t')}{|t-t'|^{\frac{\theta}2 (1-\frac2{r})}} dt dt' \leq \mathbb{E}\|f_\tau\|_{L^{q_1}}^2\\
\leq \mathbb{E}\left( \|e^{-\frac{\nu}8 (.-\tau) 2^{2j}} \textbf{1}_{\{\tau \leq \cdot\}}\|_{L^{q_2}} \|\psi\|_{L^{\bar{p}}L^{\bar{r}}} \right)^2 \leq \mathbb{E}\left(\left[\frac8{\nu q_2}\right]^{\frac{1}{q_2}} 2^{-\frac{2j}{q_2}} \|\psi\|_{L^{\bar{p}}L^{\bar{r}}} \right)^2,
\label{estimJen3}
\end{multline}
for $q_2\in [1,\infty]$ chosen so that $\frac1{q_2}+\frac1{\bar{p}}=\frac1{q_1}$, that is $\frac{1}{q_2}=\frac{1}{p}-\frac{\theta}{4}(1-\frac{2}{r})$.
\begin{rem}
 \sl{As we want $q_2\geq 1$ we need $p \leq \frac{4}{\theta (1-\frac{2}{r})}$.}
\end{rem}
Plugging \eqref{estimJen3} and \eqref{estimJen2} in \eqref{estimJen1}, we obtain that (also using that $\|\psi\|_{L^{\bar{p}}L^{\bar{r}}} \leq 1$):
\begin{multline}
 \left(\int_0^\infty  \mathbb{K}(\tau)^{\frac{\bar{k}}2} d\tau \right)^{\frac1{\bar{k}}} \leq \left(\mathbb{E}\left[\frac8{\nu q_2}\right]^{\frac{1}{q_2}} 2^{-\frac{2j}{q_2}} \right)^{1-\frac2{\bar{k}}}\\
 \times \left(\int_0^\infty \int_0^\infty \int_0^\infty \frac{\|\psi(t')\|_{L^{\bar{r}}} \|\psi(t)\|_{L^{\bar{r}}}}{|t-t'|^{\frac{\theta}2 (1-\frac2{r})}} e^{-\frac{\nu}8 (1+\frac{\bar{k}}2) 2^{2j} (t+t'-2\tau)} \textbf{1}_{\{\tau \leq \min(t,t')\}} d\tau dt dt'\right)^{\frac1{\bar{k}}}.
\end{multline}
Computing the integral in $\tau$, using the fact that $t+t'-2\min(t,t')=|t-t'|$, and introducing $g(t)=\|\psi(t)\|_{L^{\bar{r}}} \textbf{1}_{\{t\geq 0\}}$ and $W(t)=\frac{e^{-\frac{\nu}8 (1+\frac{\bar{k}}2) 2^{2j} |t|}}{|t|^{\frac{\theta}2 (1-\frac2{r})}}$, we get that:
\begin{multline}
 \left(\int_0^\infty  \mathbb{K}(\tau)^{\frac{\bar{k}}2} d\tau \right)^{\frac1{\bar{k}}} \leq \left(\mathbb{E}\left[\frac8{\nu q_2}\right]^{\frac{1}{q_2}} 2^{-\frac{2j}{q_2}} \right)^{1-\frac2{\bar{k}}} \left(\frac4{\nu (1+\frac{\bar{k}}2)} 2^{-2j}\int_\R \int_\R W(t-t') g(t) g(t') dt dt'\right)^{\frac1{\bar{k}}}\\
 \leq \left(\mathbb{E}\left[\frac8{\nu q_2}\right]^{\frac{1}{q_2}} 2^{-\frac{2j}{q_2}} \right)^{1-\frac2{\bar{k}}} \left(\frac4{\nu (1+\frac{\bar{k}}2)} 2^{-2j}\right)^{\frac1{\bar{k}}} (\|g\|_{L^{\bar{p}}} \|W*g\|_{L^p})^{\frac1{\bar{k}}}.
\end{multline}
If $p\geq 2$ then $\|W*g\|_{L^p} \leq \|W\|_{L^\frac{p}2} \|g\|_{L^{\bar{p}}}$. As soon as $\frac{p\theta}4 (1-\frac2{r}) <1$ the following integral exists and we have:
$$
 \|W\|_{L^\frac{p}2}^{\frac1{\bar{k}}}=\left(\int_\R \frac{e^{-\frac{p\nu}{16} (1+\frac{\bar{k}}2) 2^{2j} |t|}}{|t|^{\frac{p\theta}4 (1-\frac2{r})}} dt \right)^{\frac2{p\bar{k}}} =\left(\int_\R \frac{e^{-|u|}}{|u|^{\frac{p\theta}4 (1-\frac2{r})}} du \right)^{\frac2{p\bar{k}}} \left( \frac{16}{p\nu (1+\frac{\bar{k}}2)}  2^{-2j}\right)^{\frac2{p\bar{k}} -\frac{\theta}{2\bar{k}}(1-\frac2{r})},
$$
so that we end up with (using in the exponents that $\frac1{\bar{k}}=1-\frac1{k}$):
\begin{multline}
 \|\ddj f\|_{L^p L^r} \leq C_F^{\frac12-\frac1{r}} 2^{j(\frac32-\frac3{r}-\frac2{q_2}-\frac2{\bar{k}})} \ee^{\frac{\theta}4 (1-\frac2{r})} \|\ddj \Fe\|_{L^k L^2}\\
 \times \left(\mathbb{E}\left[\frac8{\nu q_2}\right]^{\frac{1}{q_2}} \right)^{1-\frac2{\bar{k}}} \left(\frac4{\nu (1+\frac{\bar{k}}2)} \right)^{\frac1{\bar{k}}} \left(\int_\R \frac{e^{-|u|}}{|u|^{\frac{p\theta}4 (1-\frac2{r})}} du \right)^{\frac2{p\bar{k}}} \left( \frac{16}{p\nu (1+\frac{\bar{k}}2)}\right)^{\frac2{p\bar{k}} -\frac{\theta}{2\bar{k}}(1-\frac2{r})}\\
 \leq \frac{C_{F, p,\theta,r,k}}{\nu^{1-\frac1{k}+\frac1{p}-\frac{\theta}4(1-\frac2{r})}} 2^{j(\frac32-\frac3{r}-\frac2{q_2}-\frac2{\bar{k}})} \ee^{\frac{\theta}4 (1-\frac2{r})} \|\ddj \Fe\|_{L^k L^2},
\end{multline}
where
\begin{multline}
 C_{F, p,\theta,r,k}=(C_F)^{\frac12-\frac1{r}} \mathbb{E}^{\frac2{k}-1}\left(\frac4{1+\frac{\bar{k}}2}\right)^{1-\frac1{k}} \left[\frac{\Big(\frac{1}{p}-\frac{\theta}{4}(1-\frac{2}{r})\Big)^{\frac2{k}-1} 2^{5-\frac2{k}}}{\left( 1+\frac{\bar{k}}2\right)^\frac2{\bar{k}} p^{2(1-\frac1{k})}}\right]^{\frac{1}{p}-\frac{\theta}{4}(1-\frac{2}{r})}\\
 \times \left(\displaystyle{\int_\R \frac{e^{-|u|}du}{|u|^{\frac{p\theta}4 (1-\frac2{r})}}}\right)^{\frac2{p}(1-\frac1{k})}
\end{multline}
Multiplying by $2^{jd}$ and summing over $j\in \Z$ ends the proof of Proposition \ref{Estimdispnu2}. $\blacksquare$

\subsubsection{Proof of Proposition \ref{estimStriPE}}
\label{preuvedispPE2}
The first line can be deduced from the second one taking $s=\frac12$ so we will focus on the last two lines.

Choosing $(d,p,r,q, \theta)=(0,\frac2{1-s}, 6, 2, 3(\frac12+\delta-s))$, thanks to Propositions \ref{injectionLr},  \ref{Propermut} (that applies when $\frac2{1-s}\geq 2$, which is true when $s\in[0,1[$) and \ref{Estimdispnu}, we have
$$
\|\Weh\|_{L_T^{\frac2{1-s}} L^6} \leq C \|\Weh\|_{L_T^{\frac2{1-s}} \dot{B}_{6,2}^0} \leq C \|\Weh\|_{\tilde{L}_T^{\frac2{1-s}} \dot{B}_{6,2}^0} \leq \frac{C}{\nu^{\frac{1-2\delta}4}} \ee^{\frac12(\frac12+\delta-s)} \|\Uoosc\|_{\dot{H}^{\frac12+\delta}}.
$$
Note that the condition $\theta\in[0,1]$ (respectively $p\in[1, \frac4{\theta(1-\frac2{r})}]$) is satisfied for any $s\in[\frac12, \frac12+\eta \delta]$ if and only if $\delta\leq \frac13$ (respectively $\delta\leq \frac12$). These conditions are true for any $s\in[\frac12-\eta\delta, \frac12[$ if and only if $\delta(1+\eta)\leq \frac13$.

With the same coefficients, and choosing $k=2$ in Proposition \ref{Estimdispnu2}, we obtain that:
$$
\|\Wei\|_{L_T^{\frac2{1-s}} L^6} \leq \frac{C}{\nu^{\frac{3-2\delta}4}} \ee^{\frac12(\frac12+\delta-s)} \|G^b\|_{L_T^2\dot{H}^{-\frac12+\delta}},
$$
and we conclude thanks to Proposition \ref{PropestimG}. The fact that $\theta\leq 1$ is satisfied with the same conditions, the condition for $p$ turns into $\delta<\frac12$. All of these are true when $\delta\leq \frac16$ and $\eta\leq 1$.\\

With the same arguments, choosing $(d,p,r,q, \theta)=(1,\frac2{\frac32-s},3,2, 6(\frac12+\delta-s))$, leads to
$$
\begin{cases}
 \|\n \Weh\|_{L_T^{\frac2{\frac32-s}}L^3} \leq \frac{C}{\nu^{\frac{1-\delta}2}} \ee^{\frac12(\frac12+\delta-s)} \|\Uoosc\|_{\dot{H}^{\frac12+\delta}},\\
 \|\n \Wei\|_{L_T^{\frac2{\frac32-s}}L^3} \leq \frac{C}{\nu^{\frac{2-\delta}2}} \ee^{\frac12(\frac12+\delta-s)} \|G^b\|_{L_T^2\dot{H}^{-\frac12+\delta}}.
\end{cases}
$$
And the result follows. Note that $\theta\in[0,1]$ now requires $\delta\leq \frac16$ and $p\in [2,\frac4{\theta(1-\frac2{r})}]$ when $\delta\leq 1$ ($\delta<1$ in the second case), so that when $\delta\leq \frac16$ all the above conditions are satisfied for $s\in[\frac12, \frac12+\eta \delta]$.
\\

To prove the last point, let us emphasize that, under the additional assumption $\Uoqg \in \dot{H}^{\frac12-\delta}$ (see $(H_5)$) we can now bound $G^b$ exactly as in \cite{FCPAA} (see (2.23)) and will not need anymore to split into $\Wei+\Weh$. For all $s\in[\frac12-\delta,\frac12+\delta]$:
\begin{multline}
  \|G^b\|_{L_T^1 \Hs } \leq C_F \|\n \tUqg\|_{L^2 \dot{H}^{\frac 12-\delta}}^{\frac12} \|\n \tUqg\|_{L^2 \dot{H}^{\frac 12+\delta}}^{\frac12} \|\n \tUqg\|_{L^2 \dot{H}^{s}}\\
  \leq \frac{\|\tUoqg\|_{\dot{H}^{\frac12-\delta} \cap \dot{H}^{\frac12+\delta}}^2}\nu e^{\frac{C}\nu \CT} \leq \frac{\Co^2}\nu e^{\frac{C}\nu \CT},
  \label{EstimGbis}
\end{multline}
so that $\Wei$ can be estimated through the same Strichartz estimates as $\Weh$. But when $s\in[\frac12-\eta\delta, \frac12]$, we have $\frac2{2-s}<2$ so we cannot bound the $L_T^{\frac2{2-s}}L^\frac3{s}$-norm with the $\tilde{L}_T^{\frac2{2-s}}\dot{B}_{\frac3{s},2}^0$-norm anymore and instead we write for $(d,p,r,q)=(1,\frac2{2-s},\frac3{s},1)$ and for any $\theta\in[0,1]$ (as $\frac2{2-s}\geq 1$ when $s\in [0,2[$):
$$
\|\n \Weh\|_{L_T^{\frac2{2-s}}L^\frac3{s}} \leq \|\n \Weh\|_{\tilde{L}_T^{\frac2{2-s}}\dot{B}_{\frac3{s},1}^0} \leq \frac{C_{F,s,\theta}}{\nu^{1-\frac{s}2-\frac{\theta}4(1-\frac{2s}3)}} \ee^{\frac{\theta}4(1-\frac{2s}3)} \|\Uoosc\|_{\dot{B}_{2,1}^{\frac12+\frac{\theta}2(1-\frac{2s}3)}}.
$$
Thanks to Lemma \ref{majBs21}, with $(\aa,\beta)=(a\frac{\theta}2(1-\frac{2s}3), b\frac{\theta}2(1-\frac{2s}3))$,
$$
\|\Uoosc\|_{\dot{B}_{2,1}^{\frac12+\frac{\theta}2(1-\frac{2s}3)}} \leq C_{a,b,\theta,s}\|\Uoosc\|_{\dot{H}^{\frac12+\frac{\theta}2(1-\frac{2s}3)(1-a)}}^{\frac{b}{a+b}} \|\Uoosc\|_{\dot{H}^{\frac12+\frac{\theta}2(1-\frac{2s}3)(1+b)}}^{\frac{a}{a+b}},
$$
and to find $a,b>0$ satisfying
\begin{equation}
 \begin{cases}
  \frac{\theta}2(1-\frac{2s}3)(1-a)= c\delta,\\
  \frac{\theta}2(1-\frac{2s}3)(1+b)= \delta,
 \end{cases}
\end{equation}
we simply set $\theta= \frac{2\delta}{1+b}\frac1{1-\frac{2s}3}$ and $a=1-(1+b)c$ for some small $b\in]0,\frac1{c}-1[$, which leads to
\begin{equation}
 \|\n \Weh\|_{L_T^{\frac2{2-s}} L^\frac3{s}} \leq \frac{C_{F,s,\delta,b,c,\Co}}{\nu^{1-\frac{s}2-\frac{\delta}{2(1+b)}}} \ee^{\frac{\delta}{2(1+b)}-\gamma}.
\end{equation}
The exponent of $\ee$ writes:
$$
\frac{\delta}{2(1+b)}-\gamma= \frac{\delta}{2(1+b)}\left(2\eta_0 -b(1-2\eta_0)\right) \underset{b\rightarrow 0}{\longrightarrow} \eta_0\delta,
$$
such that for a given $k<1$ close to 1, we can choose $b\in]0,\frac1{c}-1[$ so small that
\begin{equation}
\frac{\delta}{2(1+b)}-\gamma= k\eta_0\delta,
\end{equation}
which gives the result. To bound $\Wei$ we use the same Strichartz estimates with the same coefficients and thanks to \eqref{EstimGbis}, we obtain the rest of the estimates.
\\
To finish, let us precise that $\theta\leq 1$ for any $s\in[\frac12-\eta\delta, \frac12]$ is equivalent to the the fact that following bound is true for $s=\frac12$
$$
\frac{2\delta}{1+b} =4(\gamma+k\eta_0 \delta) \leq 1-\frac{2s}3,
$$
which is equivalent to $\delta(1-2\eta_0(1-k))\leq \frac13$ and true as we already have $\delta\leq \frac16$, $\eta\leq 2\eta_0<1$ and $k<1$.\\

Similarly, the condition on $p$ is realized when for any $s\in[\frac12-\eta\delta, \frac12]$, we have
$$
\frac{\delta}{2(1+b)} =\gamma+k\eta_0 \delta \leq 1-\frac{s}2,
$$
which is equivalent the fact that it is satisfied for $s=\frac12$, and is equivalent to asking $\delta(1-2\eta_0(1-k))\leq \frac32$, and is also true as $\delta \leq \frac16$. $\blacksquare$

\subsection{Strichartz estimates for the rotating fluids}

\subsubsection{Statement of the results}

In this section, we will provide isotropic and anisotropic strichartz estimates for System \eqref{LRF}. Let us begin with the estimates proved by Chemin, Desjardins, Gallagher and Grenier (that we present here with our notations and without external force term):
\begin{prop} (\cite{CDGG, CDGGbook})
 \sl{For any $p\in[1,\infty]$ and any $\aa>0$ there exists a constant $C$ such that for any vector field $w_0$, any $j,k \in \Z$, if $\We$ solves \eqref{LRF} with initial data $w_0$:
 $$
 \begin{cases}
  \|\ddj \We\|_{L^p L^\infty(\R^3)} \leq C 2^{j(\frac32 -\frac2{p})} \left(\ee 2^{2j}\right)^{\frac1{4p(1+\aa)}} \|\ddj w_0\|_{L^2(\R^3)},\\
  \|\ddj \ddk \We\|_{L^p L_{h,v}^{\infty,2}} \leq C 2^{j(1 -\frac2{p})} \min\left(1,\left(\ee 2^{2j}\right)^{\frac1{4p}} 2^{\frac1{2p}(j-k)} \right)\|\ddj \ddk w_0\|_{L^2(\R^3)},
 \end{cases}
 $$
 where (for $j,k\in \Z$) $\ddj=\varphi(2^{-j}D)$ and $\ddk=\varphi(2^{-k}D_3)$ are the usual homogeneous Littlewood-Paley truncation operator, and its vertical counterpart (we refer to \cite{CDGG, CDGG2, Dragos1} for details about the anisotropic Littlewood Paley theory), and where we define for $a,b\in [1,\infty]$,
 $$
 \|f\|_{L_{h,v}^{a,b}} \overset{def}{=} \big\|\|f(x_h,.)_{L^b(\R_v)}\|\big\|_{L^a(\R_h^2)}.
 $$
 }
 \label{StriCDGG}
\end{prop}
In the series of works \cite{KLT, KLT2, LT, IMT} the authors manage to improve their Strichartz estimates from \cite{IT5, IT2} thanks to the Riesz-Thorin theorem (as in \cite{Dutrifoy2}) and the Littman theorem (see references in \cite{FCcompl}: the first one allows to turn the condition "$r>4$" into "$r>2$" whereas the second allows slightly larger upper bounds for $\delta$). We also improved our Strichartz estimates from \cite{FCPAA} thanks to the same tools tools in \cite{FCcompl} and we refer to the appendix of this article for an explaination of the improvements in the rotating fluids case. We begin with the statement of the estimates we use in the proof of Theorem \ref{ThRF}.
\begin{prop}
 \sl{\begin{enumerate}
      \item For any $d\in \R$, $m\geq2$, $\theta\in[0,1]$, and $p\in[1, \frac2{\theta(1-\frac2{m})}]$, there exists a constant $C=C_{p,\theta,m}$ such that for any divergence-free vectorfield $\woe$, the solution $\We$ of \eqref{LRFb} with initial data $\woe$ satisfies:
    \begin{equation}
  \||D|^d \We\|_{\tilde{L}^p\dot{B}_{m, q}^0} \leq \frac{C_{p,\theta,m}}{\nu^{\frac{1}{p}-\frac{\theta}{2}(1-\frac2{m})}} \ee^{\frac{\theta}{2}(1-\frac2{m})} \|\woe\|_{\dot{B}_{2, q}^{\sigma_1}},
\end{equation}
with $\sigma_1= d+\frac32-\frac3{m}-\frac2{p}+\theta(1-\frac2{m})$.
      \item For any $d\in \R$, $m>2$, $\theta\in]0,1]$, $p\in[1, \frac4{\theta(1-\frac2{m})}]$ there exists a constant $C=C_{p,\theta,m}$ such that for any divergence-free vectorfield $\woe$, we have:
    \begin{equation}
  \||D|^d \We\|_{L^p L_{h,v}^{m,2}} \leq \frac{C_{p,\theta,m}}{\nu^{\frac{1}{p}-\frac{\theta}{4}(1-\frac2{m})}} \ee^{\frac{\theta}{4}(1-\frac2{m})} \|\woe\|_{\dot{B}_{2, 1}^{\sigma_2}},
\end{equation}
 with $\sigma_2= d+1-\frac2{m}-\frac2{p}+\frac{\theta}2 (1-\frac2{m})$.
     \end{enumerate}
\label{estimdispRF}
 }
\end{prop}

As a consequence we can state the following proposition, which allows to bound the terms involving $\We$ in \eqref{estimCas1RF}.

\begin{prop}
 \sl{
 \begin{enumerate}
  \item Under the previous notations, if $\delta\leq \frac13$, for any $s\in[\frac12, \frac12+\eta \delta]$, there exists a constant $C=C(\delta,s)>0$ such that:
 $$
 \begin{cases}
\displaystyle{\|\We\|_{L_T^4 L^6} +\nu^{\frac14} \|\n \We\|_{L_T^2 L^3} \leq \frac{C \Co}{\nu^{\frac{1-2\delta}4}} \ee^{\eta_0 \delta}},\\
  \displaystyle{\|\We\|_{L_T^{\frac2{1-s}}L^6} \leq\frac{C \Co}{\nu^{\frac{1-2\delta}4}} \ee^{\frac12(\frac12+2\eta_0 \delta-s)}}.
 \end{cases}
 $$
 \item If $\delta\leq \frac14$, for any $k\in]0,1[$ (as close to 1 as we wish) and $s\in[\frac12, \frac12+\eta \delta]$, there exist $C=C(\delta,s,\gamma)>0$ such that:
 $$
 \displaystyle{\|\We\|_{L_T^{\frac2{1-s}}L_{h,v}^{\infty, 2}} +\nu^{\frac12} \|\n \We\|_{L_T^{\frac2{\frac32-s}}L_{h,v}^{4, 2}} \leq \frac{C \Co}{\nu^{\frac12\left(1-s-2\gamma-k(\frac12+\delta-s)\right)}} \ee^{\frac{k}2(\frac12+2\eta_0\delta-s)}.}
 $$
 \item The previous estimates remain valid for any $s\in[\frac12-\eta \delta, \frac12+\eta \delta]$ when $\eta\leq 2\eta_0 \min(1, \frac1{k}-1)$ but the norm in $L_T^{\frac2{\frac32-s}}L_{h,v}^{4, 2}$ has to be replaced by:
 $$
 \|\n \We\|_{L_T^{\frac2{2-s}}L_{h,v}^{\infty, 2}} \leq \frac{C \Co}{\nu^{\frac12\left(2-s-2\gamma-k(\frac12+\delta-s)\right)}} \ee^{\frac{k}2(\frac12+2\eta_0\delta-s)}.
 $$
 \end{enumerate}
 }
 \label{estimStriRF}
\end{prop}

\subsubsection{Proof of Proposition \ref{estimdispRF}}

As explained in \cite{FCcompl}, the proof of the first point is globally the same as in section 2.2 from the cited article. The only difference is that, as the hessian enjoys better properties in the case of the rotating fluids, the following estimate
$$
 \|L_j(\sigma)g\|_{L^q} \leq (C_F)^{1-\frac2{q}} \frac{2^{3j(1-\frac2{q})}}{|\sigma|^{\frac{\theta}2 (1-\frac2{q})}} \|g\|_{L^{\bar{q}}},
$$
is replaced by
$$
 \|L_j(\sigma)g\|_{L^q} \leq (C_F)^{1-\frac2{q}} \frac{2^{3j(1-\frac2{q})}}{|\sigma|^{\theta (1-\frac2{q})}} \|g\|_{L^{\bar{q}}},
$$
so that Point 1 provides a similar estimates as in Proposition \ref{Estimdispnu}, but with $\theta$ replaced by $2\theta$.

Let us focus on the second point, which extends the anisotropic estimates from \cite{CDGG}. Assume that $\We$ solves \eqref{LRF} with initial data $\woe$. For any $p,m\in [1,\infty]$ and any fixed $j\in \Z$, we can write:
\begin{equation}
 \|\ddj \We\|_{L^p L_{h,v}^{m,2}} \leq \Sum_{k\leq j+1} \|\ddj \ddk \We\|_{L^p L_{h,v}^{m,2}},
\end{equation}
and, defining
$$
\cB\overset{def}{=}\{\psi \in \cC_0^\infty (\R_+\times \R^3, \R), \quad \|\psi\|_{L^{\bar{p}}(\R_+, L_{h,v}^{\bar{m},2})}\leq 1\},
$$
we have (with the same truncation function $\varphi_1$ as in Section \ref{preuvedispPE})
$$
\|\ddj \ddk \We\|_{L^p L_{h,v}^{m,2}}= C \sup_{\psi \in \cB} \int_0^\infty \int_{\R^3} e^{-\nu t|\xi|^2+i\frac{t}{\ee}\frac{\xi_3}{|\xi|}} \widehat{\ddj \ddk \woe}(\xi) \varphi_1(2^{-j}\xi) \varphi_1(2^{-k}\xi_3) \hat{\psi}(t,\xi) d\xi dt. 
$$
Following the very same steps as in \cite{FCcompl} (and Section \ref{preuvedispPE}) we get
\begin{multline}
 \|\ddj \ddk \We\|_{L^p L_{h,v}^{m,2}}\leq C \|\ddj \ddk \woe\|_{L^2}\\
 \times \sup_{\psi \in \cB} \left(\int_0^\infty \int_0^\infty \|\psi(t)\|_{L_{h,v}^{\bar{m},2}} \|\Gjk \big(\frac{t-t'}{\ee}, \nu(t+t')\big) \bar{\psi}(t')\|_{L_{h,v}^{m,2}} dt dt' \right)^{\frac12},
 \label{estimAniso1}
\end{multline}
where for any $\tau, \sigma$ and any function $g$,
$$
\Gjk(\tau, \sigma) g = \cF^{-1} \left(e^{-\sigma|\xi|^2+i \tau\frac{\xi_3}{|\xi|}} \varphi_1(2^{-j}\xi)^2 \varphi_1(2^{-k}\xi_3)^2 \hat{g}(\xi) \right).
$$
The Plancherel identity implies that (as in Section \ref{preuvedispPE}, $\varphi_1$ is supported in the annulus centered at zero and of radii $\frac12$ and 3.)
\begin{equation}
 \|\Gjk(\tau, \sigma)\|_{L^2 \rightarrow L^2} =\|\Gjk(\tau, \sigma)\|_{L_{h,v}^{2,2} \rightarrow L_{h,v}^{2,2}} \leq C e^{-\frac{\sigma}4 2^{2j}}.
 \label{InterpA}
\end{equation}
Moreover, thanks to Lemma 3 from \cite{CDGG}, we also have for any $\theta\in[0,1]$:
\begin{equation}
 \|\Gjk(\tau, \sigma)\|_{L_{h,v}^{1,2} \rightarrow L_{h,v}^{\infty,2}} \leq C \min(1, |\tau|^{-\frac12} 2^{j-k}) 2^{2j} e^{-\frac{\sigma}4 2^{2j}} \leq C |\tau|^{-\frac{\theta}2} 2^{\theta(j-k)} 2^{2j} e^{-\frac{\sigma}4 2^{2j}}.
 \label{InterpB}
\end{equation}
Gathering \eqref{InterpA} and \eqref{InterpB} and using the Riesz-Thorin theorem, we end-up for any $m\in[2,\infty]$ and $\theta \in[0,1]$ with:
\begin{equation}
 \|\Gjk(\tau, \sigma)\|_{L_{h,v}^{\bar{m},2} \rightarrow L_{h,v}^{m,2}} \leq C e^{-\frac{\sigma}4 2^{2j}} |\tau|^{-\frac{\theta}2(1-\frac2{m})} 2^{\theta(j-k)(1-\frac2{m})} 2^{2j(1-\frac2{m})}.
\end{equation}
Therefore, plugging this estimate into \eqref{estimAniso1},
\begin{multline}
 \|\ddj \ddk \We\|_{L^p L_{h,v}^{m,2}} \leq C \|\ddj \ddk \woe\|_{L^2}\\
 \times \ee^{\frac{\theta}4(1-\frac2{m})} 2^{\frac{\theta}2 (j-k)(1-\frac2{m})} 2^{j(1-\frac2{m})} \sup_{\psi \in \cB} \left(\int_0^\infty \int_0^\infty \frac{g(t)g(t')}{|t-t'|^{\frac{\theta}2(1-\frac2{m})}} dt dt' \right)^{\frac12},
\end{multline}
where we put $g(t)=\|\psi(t)\|_{L_{h,v}^{\bar{m},2}} e^{-\frac{\nu}4 t 2^{2j}}$. Following the same steps as in \cite{FCPAA, FCcompl} we end up, thanks to the Hardy-Littlewood estimate, with:
\begin{equation}
 \|\ddj \ddk \We\|_{L^p L_{h,v}^{m,2}} \leq C \|\ddj \woe\|_{L^2} \ee^{\frac{\theta}4(1-\frac2{m})} 2^{\frac12 \theta(j-k)(1-\frac2{m})} 2^{j(1-\frac2{m}-\frac2{q_2})} \left(\frac4{q_2 \nu}\right)^\frac1{q_2},
\end{equation}
with $q_2$ defined by $\frac1{q_2}=\frac1{p}-\frac{\theta}4 (1-\frac2{m})$ (the condition on $p$ comes from the fact that we ask $q_2\in[1,\infty]$). Next, summing for $k\leq j+1$ (which explains why we ask $m>2$ and $\theta>0$) we get that:
\begin{equation}
 \|\ddj \We\|_{L^p L_{h,v}^{m,2}} \leq \frac{C_{p,\theta,m}}{\nu^{\frac1{p}-\frac{\theta}4 (1-\frac2{m})}} \|\ddj \woe\|_{L^2} \ee^{\frac{\theta}4(1-\frac2{m})} 2^{j\left(1-\frac2{m}-\frac2{p}+\frac{\theta}2 (1-\frac2{m})\right)}.
\end{equation}
Multiplying by $2^{jd}$ and summing over $j\in \Z$ concludes the proof of Point 2. $\blacksquare$

\subsubsection{Proof of Proposition \ref{estimStriRF}}

Similarly to the proof of Proposition \ref{estimStriPE}, we use here Proposition \ref{estimdispRF}. Point 1 is proven choosing $(d,p,m,q, \theta)\in\{(1,2,3,2, 3\delta),(0,\frac2{1-s},6,2, \frac32(\frac12+\delta-s))\}$ and we get:
$$
 \begin{cases}
\displaystyle{\|\n \We\|_{L_T^2 L^3} \leq \frac{C}{\nu^{\frac{1-2\delta}4}} \ee^{\frac{\delta}2} \|\woe\|_{\dot{H}^{\frac12+\delta}}},\\
  \displaystyle{\|\We\|_{L_T^{\frac2{1-s}}L^6} \leq\frac{C}{\nu^{\frac{1-\delta}2}} \ee^{\frac12(\frac12+\delta-s)}}\|\woe\|_{\dot{H}^{\frac12+\delta}}.
 \end{cases}
$$
In the first case, the fact that $\theta\in [0,1]$ requires $\delta\leq \frac13$ and the condition $p\leq \frac2{\theta(1-\frac2{m})}$ requires $\delta\leq 1$. In the second case, the fact that these conditions are true for any $s\in[\frac12, \frac12+\eta \delta]$ require respectively $\delta\leq \frac23$ and $\delta \leq \frac12$. The second condition is also true for any $s\in[\frac12-\eta\delta, \frac12]$ and the first one is true for such $s$ when $\delta(1+\eta)\leq \frac23$, which is realized when $\delta \leq \frac14$ and $\eta\leq 1$.

Let us now turn to the anisotropic estimates from Point 2. As the summability index in the Besov spaces from the second point of Proposition \ref{estimdispRF} is equal to 1 (our estimates do not allow it to be equal to $2$), we have no choice but asking not only that the $\dot{H}^{\frac12+\delta}$-norm but also the $\dot{H}^{\frac12 +c\delta}$-norm of $\woe$ are bounded by $\ee^{-\gamma}$.
\\

Choosing $(d,p,m)=(0, \frac2{1-s},\infty)$ and using Lemma \ref{majBs21} with $(\aa,\beta)=(a\frac{\theta}2, b\frac{\theta}2)$ (with $a,b<0$) leads to:
$$
\|\We\|_{L_T^{\frac2{1-s}}L_{h,v}^{\infty, 2}} \leq \frac{C}{\nu^{\frac{1-s}2-\frac{\theta}4}} \ee^{\frac{\theta}4} \|\woe\|_{\dot{B}_{2,1}^{s+\frac{\theta}2}} \leq \frac{C}{\nu^{\frac{1-s}2-\frac{\theta}4}} \ee^{\frac{\theta}4} C_{a,b,\theta} \|\woe\|_{\dot{H}^{s+\frac{\theta}2(1-a)}}^\frac{b}{a+b} \|\woe\|_{\dot{H}^{s+\frac{\theta}2(1+b)}}^\frac{a}{a+b}.
$$
We recall that $c\in]0,1[$ is expected to be close to 1, and as was done in \cite{FCPAA, FCcompl}, we want to choose $a,b>0$ so small that:
$$
\begin{cases}\vspace{0.1cm}
 s+\frac{\theta}2(1-a) =\frac12+c\delta,\\
 s+\frac{\theta}2(1+b) =\frac12+\delta.
\end{cases}
$$
The rest is very close to what we did in Section \ref{preuvedispPE2}: for some $b>0$ to be fixed later, let us take $\theta=\frac2{1+b}(\frac12+\delta-s)$ then the existence of some $a\in]0,1[$ satisfying the other condition is equivalent to the fact that $b>0$ is so small that $(\frac12+c\delta-s)(1+b) <\frac12+\delta-s$ and in that case:
$$
a=1-(1+b)\frac{\frac12+c\delta-s}{\frac12+\delta-s},
$$
and we obtain:
\begin{equation}
 \|\We\|_{L_T^{\frac2{1-s}}L_{h,v}^{\infty, 2}} \leq \frac{C}{\nu^{\frac12\left(1-s-\frac1{1+b}(\frac12+\delta-s) \right)}} \ee^{\frac1{2(1+b)}(\frac12+\delta-s)} \|\woe\| _{\dot{H}^{\frac12+ c\delta}\cap \dot{H}^{\frac12 + \delta}}.
 \label{estimtheta}
\end{equation}
Thanks to Assumption $(H_2')$ from Theorem \ref{ThRF}, we can choose $b$ so small that:
\begin{equation}
\frac1{2(1+b)}(\frac12+\delta-s)-\gamma= \frac{k}2 (\frac12+2\eta_0 \delta-s),
\end{equation}
and plugging this into \eqref{estimtheta} gives the estimate. For this choice of $b$, the corresponding $\theta$ is in $[0,1]$ if and only if $2\gamma +k(\frac12+2\eta_0 \delta-s) \leq \frac12$, and the fact that it is true for any $s\in[\frac12, \frac12+\eta \delta]$ is equivalent to
$$
2\gamma +2k\eta_0 \delta \leq \frac12 \Longleftrightarrow \Big(1-2\eta_0(1-k)\Big) \delta \leq \frac12,
$$
which is true as soon as $\delta\leq\frac12$, $k<1$ and $2\eta_0<1$.
\\

Now $\theta\in[0,1]$ for any $s\in[\frac12-\eta \delta, \frac12+\eta \delta]$ if and only if $2\gamma +k(\frac12+2\eta_0 \delta-s) \leq \frac12$ is satisfied for $s=\frac12-\eta \delta$, which is equivalent to
\begin{equation}
 \Big(1-2\eta_0+k(\eta+2\eta_0)\Big) \delta \leq \frac12,
\label{Condaniso1}
 \end{equation}
which is true when $\delta\leq \frac12$ and $\eta \leq 2\eta_0(\frac1{k}-1)$.
\\

Let us turn to the "p-index" from Proposition \ref{estimdispRF}. From the equivalence:
\begin{equation}
 \frac2{1-s}\leq \frac{2(1+b)}{\frac12+\delta-s} =\frac2{2\gamma+k(\frac12+2\eta_0\delta-s)}
\Longleftrightarrow \delta(1-2\eta_0)+\frac{k}2 +2k\eta_0 +(1-k)s\leq 1,
\end{equation}
we get that $p\leq \frac2{\theta(1-\frac2{m})}$ for all $s\in[\frac12-\eta \delta, \frac12+\eta \delta]$, is equivalent to the fact it is true for $s=\frac12+\eta \delta$, that is
\begin{equation}
\delta\Big(1-(1-k)(2\eta_0-\eta)\Big) \leq \frac12,
 \label{Condaniso2}
\end{equation}
which is true when $\delta\leq \frac12$ and $\eta\leq 2\eta_0$.

The second term is treated choosing $(d,p,m)=(1,\frac2{\frac32-s},4)$ and $\theta=\frac2{1+b}(\frac12+\delta-s)$ with the same $a,b$ as in the previous lines, for wich the analogous conditions on $\theta,p$ require that \eqref{Condaniso1} and \eqref{Condaniso2} are true but for $\frac14$ instead of $\frac12$ in the right-hand-side, which explains that the final condition for all of them to be true for any $s\in[\frac12 -\eta \delta, \frac12+\eta \delta]$ is:
$$
\delta\leq \frac14 \quad \mbox{and} \quad\eta\leq 2\eta_0 \min(1, \frac1{k}-1).
$$
Finally, the last point is treated choosing $(d,p,m)=(1,\frac2{2-s},2)$ and $\theta=\frac2{1+b}(\frac12+\delta-s)$ with the very same choice for $b$ and conditions, which concludes the proof. $\blacksquare$

\textbf{Aknowledgements :} This work was supported by the ANR project INFAMIE, ANR-15-CE40-0011.

\end{document}